\documentclass[11pt, letterpaper]{article}
\usepackage[margin = 1in]{geometry}
\usepackage{graphicx,color}
\usepackage[cp850]{inputenc}
\usepackage[T1]{fontenc}
\usepackage[english]{babel}
\usepackage[f]{esvect}
\usepackage{amsmath, amsthm, amssymb}
\usepackage{rotating}
\usepackage{mathrsfs}
\usepackage{url}
\usepackage{cite}

\def\tablenotes{\bgroup\parfillskip=0pt plus 1fil
\leftskip=0pt\relax \rightskip=0pt
\vskip2pt\footnotesize}
\def\endtablenotes{\vskip1pt\egroup}

\newtheorem{theorem}{Theorem}[section]

\newtheorem{definition}[theorem]{Definition}

\newcommand{\sinc}{\mathrm{sinc}}

\allowdisplaybreaks

\begin{document}

\title{A Numerical Treatment of Energy Eigenvalues of Harmonic Oscillators Perturbed by a Rational Function}

\author{Philippe Gaudreau and Hassan Safouhi\footnote{Corresponding author: hsafouhi@ualberta.ca. \newline The corresponding author acknowledges the financial support for this research by the Natural Sciences and Engineering Research Council of Canada~(NSERC) - Grant 250223-2011.}\\
Mathematical Section, Campus Saint-Jean\\
University of Alberta\\
8406, 91 Street, Edmonton, Alberta T6C 4G9, Canada}

\date{}

\maketitle


\begin{abstract}
\hskip -0.6cm In the present contribution, we apply the double exponential Sinc-collocation method (DESCM) to the one-dimensional time independent Schr\"odinger equation for a class of rational potentials of the form $V(x) =p(x)/q(x)$. This algorithm is based on the discretization of the Hamiltonian of the Schr\"odinger equation using Sinc expansions. This discretization results in a generalized eigenvalue problem where the eigenvalues correspond to approximations of the energy values of the corresponding Hamiltonian. A systematic numerical study is conducted, beginning with test potentials with known eigenvalues and moving to rational potentials of increasing degree.
\end{abstract}

\vspace*{0.15cm}
{\bf Keywords}

Anharmonic oscillators. Rational potential. Sinc collocation methods. Exponential transformations.

\clearpage
\section{Introduction}
Over the last three decades, the study of perturbed quantum harmonic oscillators has been on the edge of thrilling and exciting research. Numerous methods have been developed to compute the energy eigenvalues of the Schr\"odinger equation ${\cal H}\psi =E\psi$, where the Hamiltonian ${\cal H}= -\frac{ {\rm d^{2}}}{{\rm d}x^{2}} + V(x)$. The case where the potential $V(x)$ is a rational function of the form $V(x) = x^2+ \frac{p(x)}{q(x)}$ has been the subject of much interest.  Historically, the potential:
\begin{equation}\label{formula: original potential}
V(x) =  x^2 + \dfrac{\lambda x^2}{1+g x^{2}},
\end{equation}
with $\lambda\in \mathbb{R}$ and $g>0$ was among the first rational potentials to manifest interest in the scientific community. The Schr\"odinger equation with such a potential is the analogue of a zero-dimensional field theory with a nonlinear Lagrangian, which is still of interest in particle physics today \cite{Mitra1978, Summers2012, Magnen2008, Rivasseau2012a}. Additionally, outside the realm of field theory, double-well potentials are among the most important potentials in quantum mechanics \cite{Jelic2012, Griffiths2004}. Double well potentials occur abundantly in the study of quantum information theory or quantum computing. Concisely, quantum information theory attempts to generalize the ideas of classical information theory to the quantum world. Recently, systems of two particles in double well potentials have been studied experimentally with ultracold atoms \cite{Anderlini2007, Trotzky2008}. In 2009, it was theoretically proposed that neutral atoms held in double well potentials could be used to create quantum logic gates to be used for quantum information processing \cite{Foot2009}. Recently, Murmann et al. demonstrated that the quantum state of two ultracold fermionic atoms in an isolated double-well potential was completely controllable \cite{Murmann2014}. They were able to control the interaction strength between these two particles, the tilt of the potential as well as their tunneling rates between the two wells. These experiments provide a starting point for quantum computation with neutral atoms. Hence, further investigations into quantum systems with multiple wells could be an asset in constructing an efficient and reliable quantum computers.

There is a rich and diverse collection of techniques and algorithms available in literature to compute the energy eigenvalues of the Schr\"odinger equation with the rational potential \eqref{formula: original potential}. Among others, the use of variational methods, perturbation series, and perturbed ladder operators methods have been readily used to approximate the energy eigenvalues of this potential \cite{Bessis1980, Bessis1983, Bessis1992, Fack1986, Stubbins1995}. To understand the innate structure of the potential \eqref{formula: original potential}, several exact representations for the energy eigenvalues of this potential were found for specific relations between the parameters $\lambda$ and $g$ \cite{Chaudhuri1983, Blecher1987, Filho1989, Filho1991, Flessas1982a, Flessas1981, Heading1982, Gallas1988, Marcilhacy1985, Roy1987, Roy1987a, Roy1988, Roy1990, Berghe1989, Varma1981, Whitehead1982}. Moreover, the Hill determinant method as well as the Hill determinant method with a variational parameter have been used extensively to numerically calculate the energy eigenvalues of \eqref{formula: original potential} \cite{Estrin1990a, Agrawal1993, Hautot1981}. Methods utilizing power series, Taylor series and Pad\'e approximants have also been used frequently in literature \cite{Fernandez1991, Hodgson1988, Fassari1996, Fassari1997, Kaushal1979, Lai1982}. More general numerical techniques have also been studied and applied to the potential \eqref{formula: original potential} such as finite differences, numerical integration, Runge-Kutta methods and transfer matrix methods \cite{Fack1987, Galicia1979, Ishikawa2007, Bhagwat1981, Simos1998, Simos1997, Ying2010}.

Despite the amount of attention the potential \eqref{formula: original potential} has received in literature, there have been some attempts to treat more general rational potentials. In \cite{DaCosta2008}, the potential:
\begin{equation}
V(x) = x^2 + \frac{2g^{m-1} x^{2m}}{1+ \alpha g x^{2}} \quad \textrm{where} \quad m \in \mathbb{Z}^{+},\; \alpha>0 \;\; \textrm{and} \;\;  g \geq 0,
\end{equation}
is examined. A perturbation series is used to investigate the spectrum of energy eigenvalues for this specific potential. A strong asymptotic condition of order $(m-1)$ is proved and this series is also shown to be summable in the sense of the Borel-Leroy method.

In \cite{Hislop1990}, the potential:
\begin{equation}
V(x) = \omega^2 x^2 + \frac{f(x^2)}{g(x^2)},
\end{equation}
where $f(x^2)$ and $g(x^2)$ are polynomials in $x^2$ such that $\dfrac{f(x^2)}{g(x^2)} = o(x^2)$ as $ x \to \infty$, is investigated. A method for obtaining quasi-exact solutions for the energy eigenvalues is outlined.

In \cite{Nanayakkara2012}, asymptotic expansions for the energy eigenvalues of the potential:
\begin{equation}\label{Nanayakkara}
V(x) = x^{2N} + \frac{\lambda x^{m_{1}}}{1+ g x^{m_{2}}},
\end{equation}
with $N,m_{1},m_{2}$ being arbitrary positive integers and real parameters $\lambda$ and $g$ are derived. These asymptotic expansions are subsequently shown to be related to the corresponding energy levels $n$. Numerical results are presented for a variety of potentials of the form in equation \eqref{Nanayakkara}. The numerical accuracy of the proposed asymptotic estimates increase as their corresponding energy levels increase.

In \cite{Scherrer1988}, a matrix-continued-fraction algorithm is introduced to calculate the energy eigenvalues of rational potentials of the form:
\begin{equation}
V(x) = x^2 + \dfrac{  \displaystyle \sum_{m=1}^{M} \lambda_{2m} x^{2m}}{ 1+ \displaystyle \sum_{m=1}^{M} g_{2m} x^{2m}} \qquad \textrm{and} \qquad V(x) = x^2 + \dfrac{ \displaystyle \sum_{m=1}^{2M} \lambda_{m} x^{m}}{ 1+ \displaystyle \sum_{m=1}^{2M} g_{m} x^{m}}.
\end{equation}
This method is based on an expansion of the eigenfunctions into complete sets. In addition to the proposed method, numerical results are given only for the nested potential in equation~\eqref{formula: original potential} as well as polynomial potentials (i.e. $g_{m} = 0 \quad \forall m$ ). The numerical results provided, demonstrate a high accuracy in comparison to previous results available in literature.

In \cite{Witwit1996a}, the energy levels of the one dimensional potentials:
\begin{equation}
V(x) = x^2 + \frac{\lambda x^N}{1+g x^{2}} \qquad \textrm{and} \qquad  V(x) = \mu x^{2I} \pm \frac{g x^{2N}}{1+g\alpha x^{2M}},
\end{equation}
are evaluated using a finite difference approach. Numerical results are presented for a variety of parameters and they are in good agreement with other results obtained in literature.

As can be seen by the numerous approaches which have been made to solve this problem, there is a lack of uniformity in its resolution. Moreover, the task of evaluating energy eigenvalues for any general rational potential has only been studied lightly.  While several of these methods yield excellent results for specific cases of rational potentials namely the potential in equation \eqref{formula: original potential}, it would be favorable to have one algorithm capable of handling any rational potential of the form $V(x) = \omega x^{2m} +  \dfrac{p(x)}{q(x)}$, which can deliver all eigenvalues efficiently to a high pre-determined accuracy. It is the purpose of this research work to present such a general method based on the double exponential Sinc-collocation method (DESCM)presented in \cite{Gaudreau2014a}. The DESCM starts by approximating the wave function as a series of weighted Sinc functions. By substituting this approximation in the Schr\"odinger equation and evaluating this expression at several collocation points, we obtain a generalized eigensystem where the generalized eigenvalues are approximations to the energy eigenvalues. As will be shown in this contribution, the more terms we include in our series of weighted Sinc functions, the better the approximation to the energy eigenvalues becomes. This method has numerous advantages over previous methods as it is insensitive to changes in the numerous parameters in addition to having a near-exponential convergence rate. Moreover, the matrices generated by the DESCM have useful symmetric properties which make numerical calculations of their eigenvalues much easier than standard non-symmetric matrices.

\section{General definitions, properties and preliminaries}\label{The double exponential Sinc collocation method}
The sinc function, analytic for all $z \in \mathbb{C}$ is defined by the following expression:
\begin{equation} \label{formula: sinc functions}
\textrm{sinc}(z) = \left\{ \begin{array}{cc} \dfrac{\sin(\pi z)}{\pi z} &\quad \textrm{for} \quad z \neq 0 \\[0.3cm]
1  &\quad \textrm{for} \quad z=0. \end{array} \right.
\end{equation}
The Sinc function $S(j,h)(x)$ for $h \in \mathbb{R}^{+}$ and $j \in \mathbb{Z}$ is defined by:
\begin{equation}\label{formula: Sinc function}
S(j,h)(x) = \sinc \left( \dfrac{x}{h}-j\right).
\end{equation}

It is possible to expand well-defined functions as series of Sinc functions. Such expansions are known as Sinc expansions or Whittaker Cardinal expansions.
\begin{definition}\cite{Stenger1981}
Given any function $v(x)$ defined everywhere on the real line and any $h>0$, the Sinc expansion of $v(x)$ is defined by the following series:
\begin{equation}\label{formula: Sinc expansion}
C(v,h)(x) = \sum_{j=-\infty}^{\infty} v_{j} S(j,h)(x),
\end{equation}
where $v_{j} = v(jh)$.
The symmetric truncated Sinc expansion of the function $v(x)$ is defined by the following series:
\begin{equation}\label{formula: Truncated Sinc expansion}
C_{N}(v,h)(x) = \sum_{j=-N}^{N} v_{j} \, S(j,h)(x) \qquad \textrm{for} \qquad N \in \mathbb{N}.
\end{equation}
\end{definition}
The Sinc functions defined in \eqref{formula: Sinc function} form an interpolatory set of functions with the discrete orthogonality property:
\begin{equation}
S(j,h)(kh)  =  \delta_{j,k}   \qquad \textrm{for} \qquad j,k \in \mathbb{Z},
\end{equation}
where $\delta_{j,k}$ is the Kronecker delta function. In other words, $C_{N}(v,h)(x) = v(x)$ at all the Sinc grid points given by $x_{k} = kh$.

A major motivation behind the use of Sinc expansions stems from numerical integration. As it happens, integration of Sinc expansions have a direct connection with the composite trapezoidal rule. It is well known from the Euler-Maclaurin summation formula that the error in approximating the integral of a measurable function $v$ defined on the interval $(a,b)$ by the composite trapezoidal rule is given by:
\begin{equation}\label{Trapezoid error}
h\sum_{k=1}^{N-1} \dfrac{v(x_{k-1})+v(x_{k})}{2} - \int_{a}^{b} v(x) {\rm d}x \sim \sum_{l=1}^{\infty} h^{2l} \dfrac{B_{2l}}{(2l)!} \left( v^{2l-1}(b) - v^{2l-1}(a)  \right),
\end{equation}
where $B_{2l}$ are the Bernoulli numbers and $N$ is the number of points in the interval of integration, from $x_0=a$ to $x_N=b$. Hence, for a periodic function or a function for which $v^{(n)}(x) \to 0$ at the endpoints, the convergence is faster than any power of $h$. This observation has led to much research on the use of conformal maps that decay to zero at infinities. Indeed, given a conformal map $\phi^{-1}(x)$ of a simply connected domain in the complex plane with boundary points $a\neq b$ such as $\phi^{-1}(a)=-\infty$ and $\phi^{-1}(b)=\infty$ and a Lebesgue measurable function $v$ defined on the interval $(a,b)$ we can derive the following elegant result using Sinc expansions:
\begin{align}
\int_{a}^{b} v(x) {\rm d}x & =  \int_{-\infty}^{\infty} (v\circ\phi)(y)\phi^{\prime}(y) {\rm d}y \nonumber \\
& \approx \int_{-\infty}^{\infty} \sum_{j=-N}^{N} \tilde{v}_{j} \, S(j,h)(y) {\rm d}y \nonumber\\
& =  \sum_{j=-N}^{N}  \tilde{v}_{j} \, \int_{-\infty}^{\infty} S(j,h)(y) {\rm d}y \nonumber \\
& =  h \sum_{j=-N}^{N}  \tilde{v}_{j} \label{formula: Sinc: Trapezoidal rule},
\end{align}
where $\tilde{v}_{j} = (v\circ\phi)(jh)\phi^{\prime}(jh)$. As we can see, equation \eqref{formula: Sinc: Trapezoidal rule} is exactly the trapezoidal rule.

In~\cite{Stenger1981}, Stenger introduced a class of functions which are successfully approximated by a Sinc expansion. We present the definition for this class of functions bellow.
\begin{definition} \cite{Stenger1981} \label{defintion: Bd function space}
Let $d>0$ and let $\mathscr{D}_{d}$ denote the strip of width $2d$ about the real axis:
\begin{equation}
\mathscr{D}_{d} = \{ z \in \mathbb{C} : |\,\Im (z)|<d \}.
\end{equation}
For $\epsilon \in(0,1)$, let $\mathscr{D}_{d}(\epsilon)$ denote the rectangle in the complex plane:
\begin{equation}
\mathscr{D}_{d}(\epsilon) = \{z \in \mathbb{C} : |\,\Re(z)|<1/\epsilon, \, |\,\Im (z)|<d(1-\epsilon) \}.
\end{equation}
Let ${\bf B}_{2}(\mathscr{D}_{d})$ denote the family of all functions $g$ that are analytic in $\mathscr{D}_{d}$, such that:
\begin{equation}\label{formula: integral imaginary}
\displaystyle \int_{-d}^{d} | \,g(x+iy)| \textrm{d}y \to 0 \; \textrm{as} \; x \to \pm \infty \quad \textrm{and} \quad  \mathcal{N}_{2}(g,\mathscr{D}_{d}) = \displaystyle \lim_{\epsilon \to 0} \left(  \int_{\partial \mathscr{D}_{d}(\epsilon)}  |g(z)|^{2} |\textrm{d}z| \right)^{1/2} <\infty.
\end{equation}
\end{definition}

The time independent Schr\"{o}dinger equation is given by the following expression:
\begin{equation}\label{formula:Schrodinger equation}
{\cal H} \, \psi(x) \, \,=\,  E \, \psi(x),
\end{equation}
where the Hamiltonian is given by the following linear operator:
\begin{equation}
{\cal H} = -\dfrac{{\rm d}^2}{{\rm d} x^2} +V(x),
\end{equation}
where $V(x)$ is the potential energy function.
The time independent Schr\"{o}dinger equation~\eqref{formula:Schrodinger equation} can be written as the following boundary value problem:
\begin{equation}\label{formula: Schrodinger equation sturm liouville}
- \psi^{\prime \prime}(x) + V(x)\psi(x) = E \psi(x)  \qquad \textrm{with} \qquad \displaystyle \lim_{|x| \to \infty} \psi(x) = 0.
\end{equation}

Equation~\eqref{formula: Schrodinger equation sturm liouville} is similar to the Schr\"odinger equation with an anharmonic oscillator potential to which we applied successfully the DESCM~\cite{Gaudreau2015b}.

Eggert et al.~\cite{Eggert1987} demonstrate that applying an appropriate substitution to the boundary value problem~\eqref{formula: Schrodinger equation sturm liouville}, results in a symmetric discretized system when using Sinc expansion approximations. The change of variable they propose is given by:
\begin{equation}
v(x) = \left(\sqrt{ (\phi^{-1})^{\prime} } \, \psi \right) \circ \phi(x) \qquad \Longrightarrow \qquad  \psi(x)  =  \dfrac{ v \circ \phi^{-1}(x)}{\sqrt{ (\phi^{-1}(x))^{\prime}}},
\label{formula: EggertSub}
\end{equation}
where $\phi^{-1}(x)$ a conformal map of a simply connected domain in the complex plane with boundary points $a\neq b$ such that $\phi^{-1}(a)=-\infty$ and $\phi^{-1}(b)=\infty$.

Applying the substitution~\eqref{formula: EggertSub} to~\eqref{formula: Schrodinger equation sturm liouville}, we obtain:
\begin{align}\label{formula: transformed Schrodinger equation}
\hat{\mathcal{H}} \, v(x) & =  - v^{\prime \prime}(x) + \tilde{V}(x) v(x)
%\nonumber\\ &
\, = \, E (\phi^{\prime}(x))^{2} v(x) \quad \textrm{with} \quad \lim_{ |x| \to \infty} v(x) = 0,
\end{align}
where:
\begin{equation}
\tilde{V}(x)  =  - \sqrt{\phi^{\prime}(x)} \, \dfrac{{\rm d}}{{\rm d} x} \left( \dfrac{1}{\phi^{\prime}(x)} \dfrac{{\rm d}}{{\rm d} x}( \sqrt{\phi^{\prime}(x)})  \right) + (\phi^{\prime}(x))^{2} V(\phi(x)).
\end{equation}

\section{The development of the method}

A function $\omega(x)$ is said to decay double exponentially at infinities if there exist positive constants $A, B,\gamma$ such that:
\begin{equation}
|\, \omega(x)| \leq A \exp( -B \exp(\gamma|\,x|)) \qquad \textrm{for} \qquad x\in \mathbb{R}.
\end{equation}

The double exponential transformation is a conformal mapping $\phi(x)$ which allows for the solution of~\eqref{formula: transformed Schrodinger equation} to have double exponential decay at both infinities.

In order to implement the DESCM, we start by approximating the solution of~\eqref{formula: transformed Schrodinger equation} by a truncated Sinc expansion~\eqref{formula: Truncated Sinc expansion}.
Inserting~\eqref{formula: Truncated Sinc expansion} into~\eqref{formula: transformed Schrodinger equation}, we obtain the following system of equations:
\begin{align}\label{formula: index solution}
\hat{\mathcal{H}} \, C_{N}(v,h)(x_{k}) & = \displaystyle \sum_{j=-N}^{N} \left[ - \dfrac{{\rm d}^2}{{\rm d} x_{k}^{2}} S(j,h)(x_{k})  + \tilde{V}(x_{k}) S(j,h)(x_{k}) \right] v_{j} \\
& = \mathcal{E} \displaystyle \sum_{j=-N}^{N} S(j,h)(x_{k}) (\phi^{\prime}(x_{k}))^{2}  v_{j} \qquad \textrm{for} \qquad k = -N, \ldots, N,
\end{align}
where the collocation points $x_{k}=kh$ and $\mathcal{E}$ is an approximation of the eigenvalue $E$ in~\eqref{formula: transformed Schrodinger equation}.

The above equation can be re-written as follows:
\begin{eqnarray}\label{formula: index solution simplified}
\hat{\mathcal{H}} \, C_{N}(v,h)(x_{k}) & = &  \sum_{j=-N}^{N} \left[ -\dfrac{1}{h^{2}} \, \delta^{(2)}_{j,k} + \tilde{V}(kh) \, \delta^{(0)}_{j,k}\right] v_{j}
 \nonumber\\ & = & \mathcal{E}  \displaystyle \sum_{j=-N}^{N} \delta^{(0)}_{j,k} (\phi^{\prime}(kh))^{2}  v_{j} \qquad {\rm for} \quad k = -N, \ldots, N,
\end{eqnarray}
where $\delta^{(l)}_{j,k}$ are given by~\cite{Stenger1979}:
\begin{equation}
\delta^{(l)}_{j,k} =  h^{l} \left. \left( \dfrac{d}{d x} \right)^{l} S(j,h)(x) \right|_{x=kh}.
\end{equation}

Equation~\eqref{formula: index solution simplified} can be represented in a matrix form as follows:
\begin{align}\label{formula: matrix solution}
\hat{\mathcal{H}} \, {\bf C}_{N}(v,h) &  = {\bf H}{\bf v} \,=\, \mathcal{E} {\bf D}^{2}{\bf v} \quad \Longrightarrow \quad ({\bf H} - \mathcal{E} {\bf D}^{2}){\bf v} \,=\, 0,
\end{align}
where:
\begin{align*}
{\bf v} & = (v_{-N},\ldots, v_{N})^{T} \qquad \textrm{and} \qquad {\bf C}_{N}(v,h)  = (C_{N}(v,h)(-Nh), \ldots, C_{N}(v,h)(Nh) )^{T}.
\end{align*}

${\bf H}$ is a $(2N+1) \times (2N+1)$ matrix with entries $H_{j,k}$ given by:
\begin{equation}\label{formula: H components}
H_{j,k} =   -\dfrac{1}{h^{2}} \, \delta^{(2)}_{j,k} + \tilde{V}(kh) \, \delta^{(0)}_{j,k} \qquad {\rm with} \quad -N \leq j,k \leq N,
\end{equation}
and ${\bf D}^{2}$ is a $(2N+1) \times (2N+1)$ diagonal matrix with entries $D^{2}_{j,k}$ given by :
\begin{equation}\label{formula: D components}
D^{2}_{j,k} =  (\phi^{\prime}(kh))^{2}  \, \delta^{(0)}_{j,k}  \qquad {\rm with} \quad -N \leq j,k \leq N.
\end{equation}

To obtain nontrivial solutions for~\eqref{formula: matrix solution}, we have to set:
\begin{equation}
\det({\bf H}-\mathcal{E}{\bf D}^{2}) = 0.
\label{formula: matrix solution-DET}
\end{equation}

In order to find an approximation of the eigenvalues of equation \eqref{formula: transformed Schrodinger equation}, we will solve the above generalized eigenvalue problem. The matrix ${\bf D}^{2}$ is diagonal and positive definite. The matrix ${\bf H}$ is the sum of a diagonal matrix and a symmetric Toeplitz matrix.  If there exits a constant $\delta>0$ such that $\tilde{V}(x)\geq \delta^{-1}$, then the matrix ${\bf H}$ is also positive definite.

In \cite[Theorem 3.2]{Gaudreau2014a}, we present the convergence analysis of DESCM which we state here in the case of the transformed Schr\"odinger equation~\eqref{formula: transformed Schrodinger equation}. The proof of the Theorem is given in~\cite{Gaudreau2014a}.
\begin{theorem}\cite[Theorem 3.2]{Gaudreau2014a} \label{theorem: convergence of eigenvalues}
Let $E$ and $v(x)$ be an eigenpair of the transformed Schr\"odinger equation:
\begin{equation}\label{transformed Sturm Liouville}
 - v^{\prime \prime}(x) + \tilde{V}(x) v(x) = \, E (\phi^{\prime}(x))^{2} v(x),
\end{equation}
where:
\begin{equation}
\tilde{V}(x)  =  - \sqrt{\phi^{\prime}(x)} \, \dfrac{{\rm d}}{{\rm d} x} \left( \dfrac{1}{\phi^{\prime}(x)} \dfrac{{\rm d}}{{\rm d} x}( \sqrt{\phi^{\prime}(x)})  \right) + (\phi^{\prime}(x))^{2} V(\phi(x)) \quad \textrm{and} \quad \lim_{ |x| \to \infty} v(x) = 0.
\end{equation}
Assume there exist positive constants $A,B,\gamma$ such that:
\begin{equation}\label{formula: xi growth condtion}
 |v(x)| \leq A \exp( - B \exp(\gamma |x|)) \qquad \textrm{for all} \qquad x \in \mathbb{R},
\end{equation}
and that $v \in {\bf B}_{2}(\mathscr{D}_{d})$ with $d \leq \dfrac{\pi}{2 \gamma}$.

Moreover, if there is a constant $\delta>0$ such that $\tilde{V}(x)\geq \delta^{-1}$ and the selection of the optimal mesh size $h$ is such that:
\begin{equation}\label{optimal mesh size}
h = \dfrac{W(\pi d \gamma N / B)}{\gamma N},
\end{equation}
where $W(x)$ is the Lambert W function.

Then, we have:
\begin{equation}\label{convergence result}
|\mathcal{E}-E|   \leq \vartheta_{v,d} \sqrt{\delta E} \left(\dfrac{N^{5/2}}{\log(N)^{2}} \right)  \exp \left(-  \dfrac{\pi d \gamma N}{\log(\pi d \gamma N/B)} \right) \quad \textrm{as} \quad N\to \infty,
\end{equation}
where $\mathcal{E}$ is the eigenvalue of the generalized eigenvalue problem and $\vartheta_{v,d}$ is a constant that depends on $v$ and $d$.
\end{theorem}
As we can see from the results obtained in Theorem~\ref{theorem: convergence of eigenvalues}, $|\mathcal{E}-E| \to 0 $ as $N\to \infty$ for all energy eigenvalues~$E$. Moreover, it is important to notice that the convergence rate of the DESCM is dependent on the strip width of analyticity $d$.

\section{The Rational Potentials} \label{Section_4}
In this section, we will apply the DESCM to the time independent Schr\"odinger equation with a rational potential. The time independent Schr\"odinger equation is defined by the following equation:
\begin{equation}
\begin{array}{ccl}
\left[ -\dfrac{{\rm d}^2}{{\rm d} x^2} +V(x)\right] \, \psi(x) \, & = & \, E \, \psi(x),
\end{array}
\end{equation}
where the potential $V(x)$ is given by:
\begin{equation} \label{formula: potential}
V(x) = \omega x^{2m} + \dfrac{p(x)}{q(x)} \quad \textrm{with} \quad (\omega,m) \in \mathbb{R}_{>0} \times \mathbb{N},
\end{equation}
where the functions $p(x)$ and $q(x)$ are polynomials in $x$. More specifically, we will investigate a subset $\mathcal{Q}$ of the set $\mathbb{R}[x]$ of rational functions of $x$ with real coefficients. The subset $\mathcal{Q}$ is defined by the following function space:
\begin{equation}\label{formula: function space V}
\mathcal{Q} = \left \{  \dfrac{p}{q} \in \mathbb{R}[x]: q(x)\,\neq\, 0 \;\; \textrm{for all} \;\; x \in \mathbb{R}, \;\;  q(0)  = 1 \;\;  \textrm{and} \;\;  \dfrac{p(x)}{q(x)}   =   o(x^{2m}) \;\;  \textrm{as} \;\; |x| \to \infty\right \}.
\end{equation}

The first condition, $q(x) \neq 0$ for all $x \in \mathbb{R}$ is set in order for $V(x)$ to be continuous for all $x\in \mathbb{R}$. However, this conditions automatically imposes that $\deg(q(x))$ be even since any odd degree polynomial has at least one real root. The second condition imposes a uniqueness on the potentials. For example, we can create the rational function:
\begin{equation}
\dfrac{p_{1}(x)}{q_{1}(x)} = \dfrac{2 + 2x + 2 x^4}{ 2 + 2 x^2} = \dfrac{1 + x + x^4}{1 +  x^2} = \dfrac{p_{2}(x)}{q_{2}(x)}.
\end{equation}
Although $p_{1}(x) \neq p_{2}(x)$ and $q_{1}(x) \neq q_{2}(x)$, we have $p_{1}(x)/q_{1}(x)  = p_{2}(x)/q_{2}(x)$. The condition $q(0) =  1$ removes this ambiguity when defining new potentials. The third condition specifies that for large $x$, the rational potential $V(x)$ exhibits behaviour similar to an anharmonic potential of the form $\omega x^{2m}$. Consequently, the potential $V(x)$ we will be investigating has the following general form:
\begin{equation}\label{formula: rational oscillator}
V(x) = \omega x^{2m} + \dfrac{\displaystyle \sum_{i=0}^{k} \lambda_{i} x^{i}   }{1 + \displaystyle \sum_{i=1}^{2l} g_{i} x^{i}} \qquad \textrm{with} \quad  k-2l<2m,
\end{equation}
It is worth noting that this current form also encompasses even potentials depending on the values of the coefficients $\{\lambda_{i}\}_{i=0}^{k}$ and $\{g_{i}\}_{i=1}^{2l}$.  To implement the DE transformation, we choose a function $\phi$ which would result in the solution of~\eqref{formula: transformed Schrodinger equation} to decay double exponentially.  Given the nature of the potential \eqref{formula: rational oscillator}, we know that the wave function will be analytic everywhere in the complex plane except for the points $z \in \mathbb{C}$ where $q(z)=0$. Since the anharmonic potential in equation \eqref{formula: rational oscillator} is analytic in $\mathbb{R}$ and grows to infinity as $x\to\pm\infty$, the wave function is also analytic and normalizable over~$\mathbb{R}$. To find the decay rate of the solution to equation ~\eqref{formula: Schrodinger equation sturm liouville} with the rational potential, we can use a similar WKB analysis as the one presented in \cite{Gaudreau2015b}. In \cite{Gaudreau2015b}, we considered the anharmonic potential:
\begin{equation}\label{formula: anharmonic potential}
V(x) = \sum_{i=1}^{m} c_{i} x^{2i} = {\mathcal O}( x^{2m} ) \quad \textrm{as} \quad |x| \to \infty.
\end{equation}
As we can see, the anharmonic potential \eqref{formula: anharmonic potential} has the same asymptotic growth as the rational potential \eqref{formula: rational oscillator}. Consequently, the function $\psi(x)$ has the same decay rate at both infinities:
\begin{equation}\label{formula: wave function asymptotic}
\psi(x) = {\cal O} \left( |x|^{-m/2} \exp \left( - \dfrac{ \sqrt{\omega}\, |x|^{m+1}}{m+1} \right) \right) \qquad \textrm{as} \qquad |x| \to \infty.
\end{equation}
Away from both infinities, the wave function $\psi(x)$ undergoes oscillatory behaviour. From equation \eqref{formula: wave function asymptotic}, we can see that the wave function $\psi(x)$ decays only single exponentially at infinities.

A possible choice for the conformal mapping is $\phi(x) = \sinh(x)$ since:
\begin{eqnarray}\label{formula: wave function asymptotic 4}
|v(x)| & = & \left| \dfrac{\psi \circ \phi(x)}{\sqrt{\phi^{\prime}(x)}} \right| \nonumber\\
& \leq & A |\sinh(x)|^{-m/2}|\cosh(x)|^{-1/2}  \exp \left( - \dfrac{ \sqrt{\omega} |\sinh(x)|^{m+1}}{m+1} \right) \nonumber\\
& \leq &  A \exp \left(  - \dfrac{ \sqrt{\omega}}{(m+1)2^{m+1}} \exp((m+1)|x|) \right),
\end{eqnarray}
for some positive constant $A$. However, as previously mentioned, Theorem \ref{theorem: convergence of eigenvalues} demonstrates that the convergence of the DESCM is dependent on the strip of width $2d\leq\dfrac{\pi}{\gamma}$ for which the function $v(x)$ is analytic. Since, we know that our original wave function $\psi(x)$ is analytic everywhere in the complex plane except for the points $z \in \mathbb{C}$ where $q(z)=0$; we wish to find a conformal mapping which will result in the solution of the transformed wave equation, $v(x)$, in \eqref{formula: transformed Schrodinger equation} to be analytic up to the maximal strip width $d=\dfrac{\pi}{2\gamma}$. To choose this optimal mapping $\phi$, we will utilize a result presented in \cite{Slevinsky2014a}.  Slevinsky et al. investigate the use of conformal maps for the acceleration of convergence of the trapezoidal rule and Sinc numerical methods \cite{Slevinsky2014a}. The conformal map they propose is constructed as a polynomial adjustment to a sinh map, and allows for the treatment of a finite number of singularities in the complex plane. The polynomial adjustments achieve this goal by locating singularity pre-images on the boundary of the widest allowable strip $\partial \mathscr{D}_{\frac{\pi}{2\gamma}}$. The map they propose has the form:
\begin{equation}\label{formula:optimal-map}
\phi(x)= u_{0}\sinh(x) + \sum_{j=1}^{n} u_{j} x^{j-1} \qquad \textrm{with} \quad u_{0} > 0,
\end{equation}
for the $(n+1)$ coefficients $\{u_{k} \}_{k=0}^{n}$ to be determined given a finite set of singularities $\{\delta_{k} \pm \epsilon_{k}\}_{k=1}^{n}$. Equation \eqref{formula:optimal-map} still grows single-exponentially. Hence, the transformed wave equation in \eqref{formula: transformed Schrodinger equation} will still result in a double exponential variable transformation. Indeed, using this conformal mapping, we see that for all $x \in \mathbb{R}$:
\begin{equation}\label{formula: wave function asymptotic 2}
|v(x)|  \leq  A\exp \left(  - \dfrac{ \sqrt{\omega}}{m+1} \left( \dfrac{u_{0}}{2} \right)^{m+1} \exp((m+1)|x|)\right),
\end{equation}
for some positive constant $A$. From equation \eqref{formula: wave function asymptotic 2}, we deduce that:
 \begin{align}
\gamma  = m+1 \qquad \textrm{and} \qquad  B  = \dfrac{ \sqrt{\omega}}{m+1} \left( \dfrac{u_{0}}{2} \right)^{m+1}.
\label{formula: B optimal}
\end{align}

Concisely, the algorithm presented in \cite{Slevinsky2014a} is as follows; given a finite set of singularities $\{\delta_{k} \pm \epsilon_{k}\}_{k=1}^{n}$, we wish to solve the following system of complex equations:
\begin{equation}\label{formula: system of equations}
\phi\left(x_{k} + i\dfrac{\pi}{2\gamma} \right) = \delta_{k} + i \epsilon_{k} \qquad \textrm{for} \quad k = 1,\ldots,n.
\end{equation}

This is a system of $n$ complex equations for the $2n+1$ unknowns $\{u_{k} \}_{k=0}^{n}$ and the real parts of the pre-images of the singularities $\{x_{k} \}_{k=1}^{n}$. Since there is one more unknown than equations, we can maximize the value of $u_{0}$ which is proportional to $B$ in equation \eqref{formula: B optimal}. By summing all $n$ equations in \eqref{formula: system of equations} and solving for $u_{0}$, we obtain the following non-linear program:
\begin{align}\label{formula: nonlinear program}
\textrm{maximize} \,\, u_{0} \,= & \left( \dfrac{  \displaystyle \sum_{k=1}^{n} \left[ \epsilon_{k} - \Im \left\{ \sum_{j=1}^{n} u_{j} \left( x_{k}+ i \dfrac{\pi}{2\gamma} \right)^{j-1} \right\} \right]}{\sin\left( \dfrac{\pi}{2\gamma}\right)  \displaystyle \sum_{k=1}^{n}\cosh(x_{k}) } \right) \\
\textrm{subject to}   &\quad \phi\left(x_{k} + i\dfrac{\pi}{2\gamma} \right) = \delta_{k} + i \epsilon_{k} \qquad \textrm{for} \quad k = 1,\ldots,n, \nonumber
\end{align}
where $\Im \left\{z\right\}$ stands for the imaginary part of $z$.

For more information of the implementation of this algorithm, we refer the readers to \cite{Slevinsky2014a}. Implementing this algorithm guarantees that all complex singularities of the potential in \eqref{formula: rational oscillator} will lie on the lines $\pm i \dfrac{\pi}{2\gamma}$ implying that our transformed solution $v(x)\in {\bf B}_{2} (\mathscr{D}_{\frac{\pi}{2\gamma}})$. Hence, the convergence rate of the DESCM will now be given by:
\begin{equation}\label{Optimized convergence result}
|\mathcal{E}-E|   \leq \vartheta_{v} \sqrt{\delta E} \left(\dfrac{N^{5/2}}{\log(N)^{2}} \right)  \exp \left(-  \dfrac{\pi^2 N}{2\log(\pi^2 N/2B)} \right) \quad \textrm{as} \quad N\to \infty.
\end{equation}

\section{Numerical discussion}
In this section, we present numerical results for the energy values for rational potentials discussed in the previous section. All calculations are performed using the programming language Julia~\cite{Bezanson2012} in double precision. The eigenvalue solvers in Julia utilize the linear algebra package {\it LAPACK}~\cite{Anderson1999}. In order to obtain the optimal conformal map in equation~\eqref{formula:optimal-map}, the non-linear program in~\eqref{formula: nonlinear program} is solved using the package {\it DEQuadrature} developed by Slevinksy and avbailable at {\it https://github.com/MikaelSlevinsky/DEQuadrature.jl}. The matrices ${\bf H}$ and ${\bf D}^{2}$ are constructed using equations~\eqref{formula: H components} and \eqref{formula: D components}. To produce our figures, we use the Julia package {\it PyPlot}.

In \cite{Chaudhuri1983, Blecher1987, Filho1989, Filho1991, Flessas1982a, Flessas1981, Heading1982, Gallas1988, Marcilhacy1985, Roy1987, Roy1987a, Roy1988, Roy1990, Berghe1989, Varma1981, Whitehead1982}, several authors presented exact solutions for potentials of the form:
\begin{equation}\label{formula: exact potential}
V(x) = x^2 + \dfrac{\lambda(g)x^2}{1+gx^2},
\end{equation}
where $\lambda$ was dependant on $g$. A few examples are presented bellow:
\begin{equation}\label{formula: true value energy}
\begin{array}{lllll}
\lambda_{1}(g) & = & -2g(2+g)  & \Rightarrow &  E_{0} = 5+ \lambda_{1}(g)/g \\
\lambda_{2}(g) & = & -2g(2+3g) & \Rightarrow &  E_{1} = 7+\lambda_{2}(g)/g \\
\lambda_{3}(g) & = & -g\left(7g+6-\sqrt{25g^2-12g+4}\right) &\Rightarrow &  E_{2} = 9+ \lambda_{3}(g)/g \\
\lambda_{4}(g) & = & -g\left(13g+6-\sqrt{49g^2-4g+4}\right) &\Rightarrow &  E_{3} = 11+\lambda_{4}(g)/g.
\end{array}
\end{equation}

\begin{figure}[!h]
\begin{center}
\begin{tabular}{cccc} \includegraphics[width=0.40\textwidth]{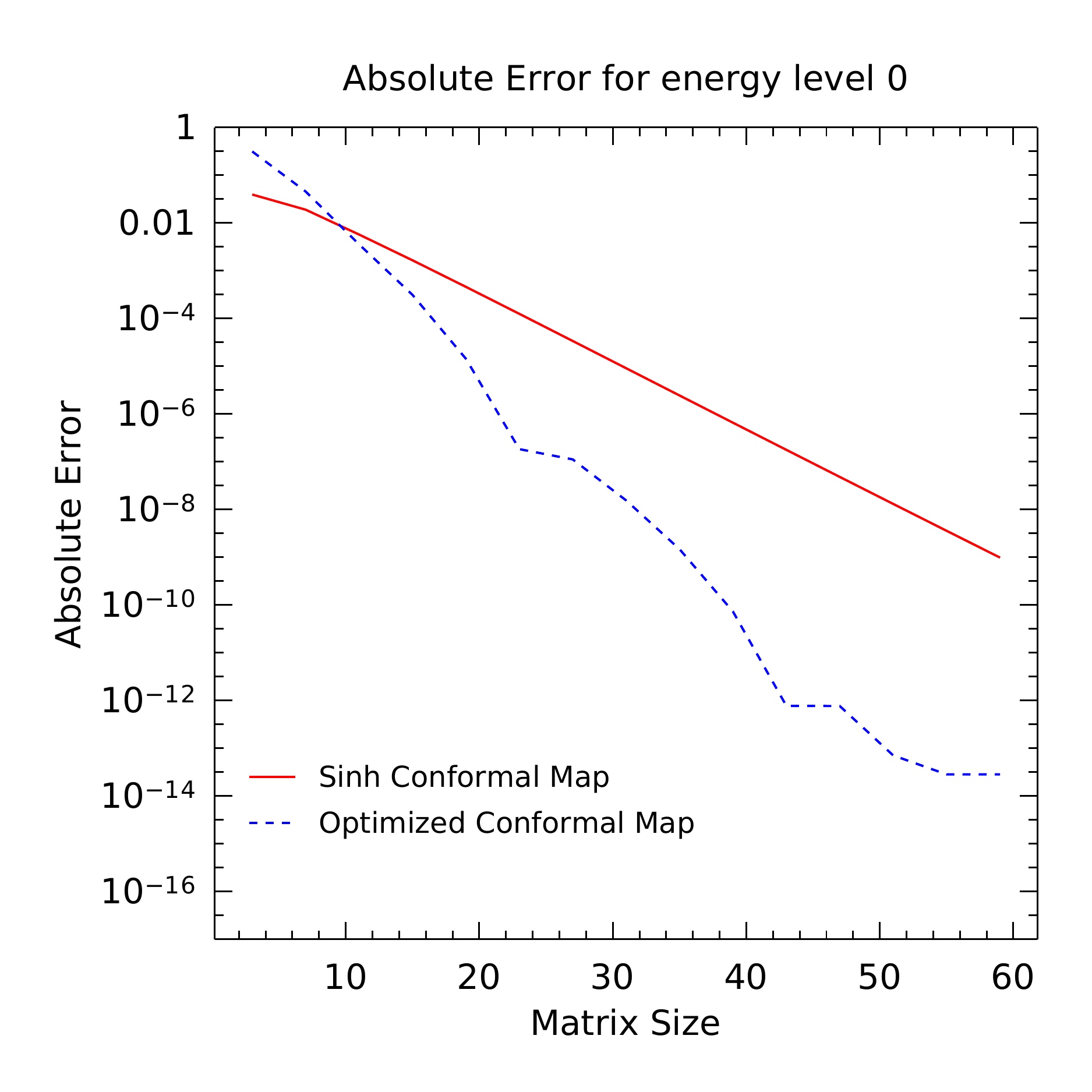} &  \includegraphics[width=0.40\textwidth]{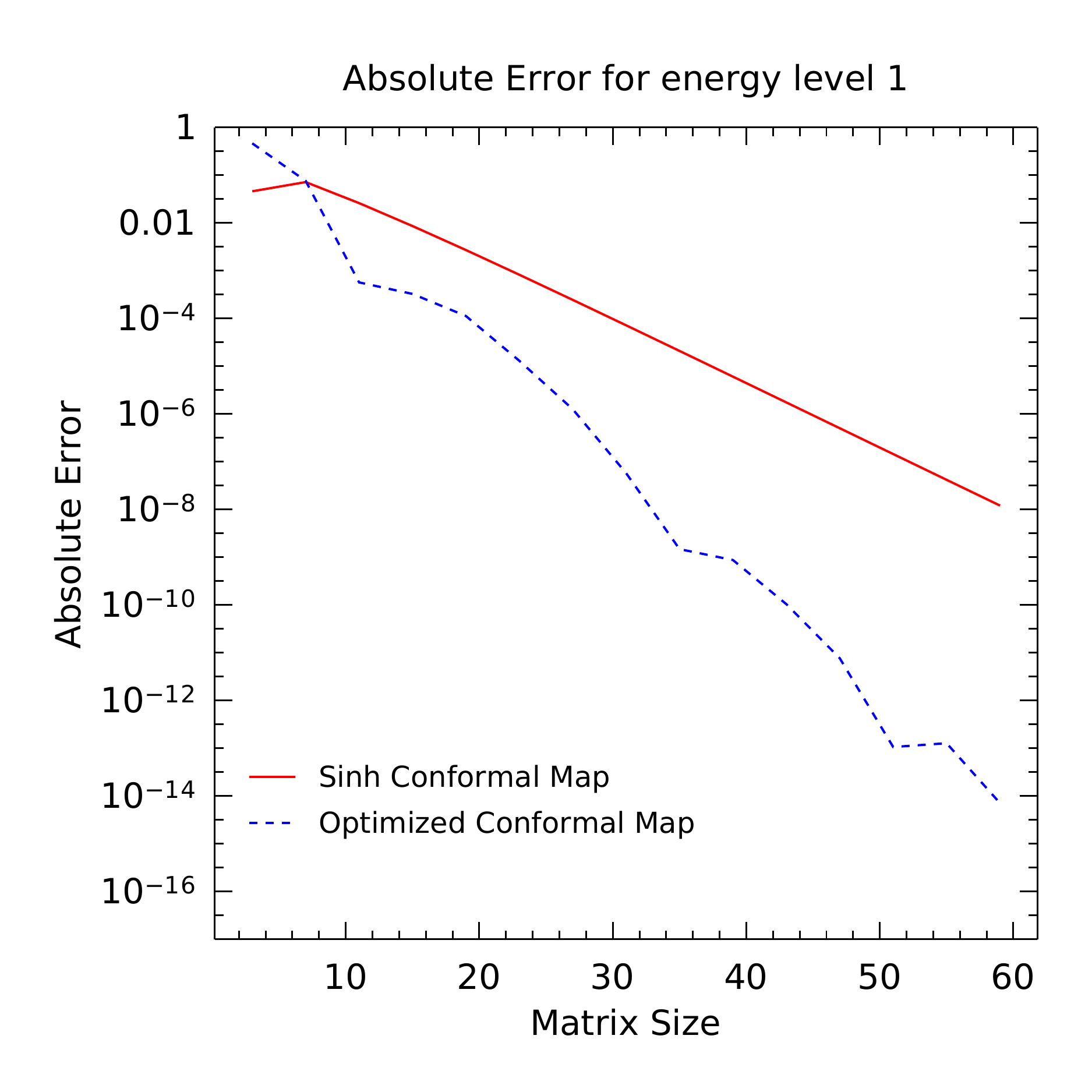} \\
(a) & (b) \\
\includegraphics[width=0.40\textwidth]{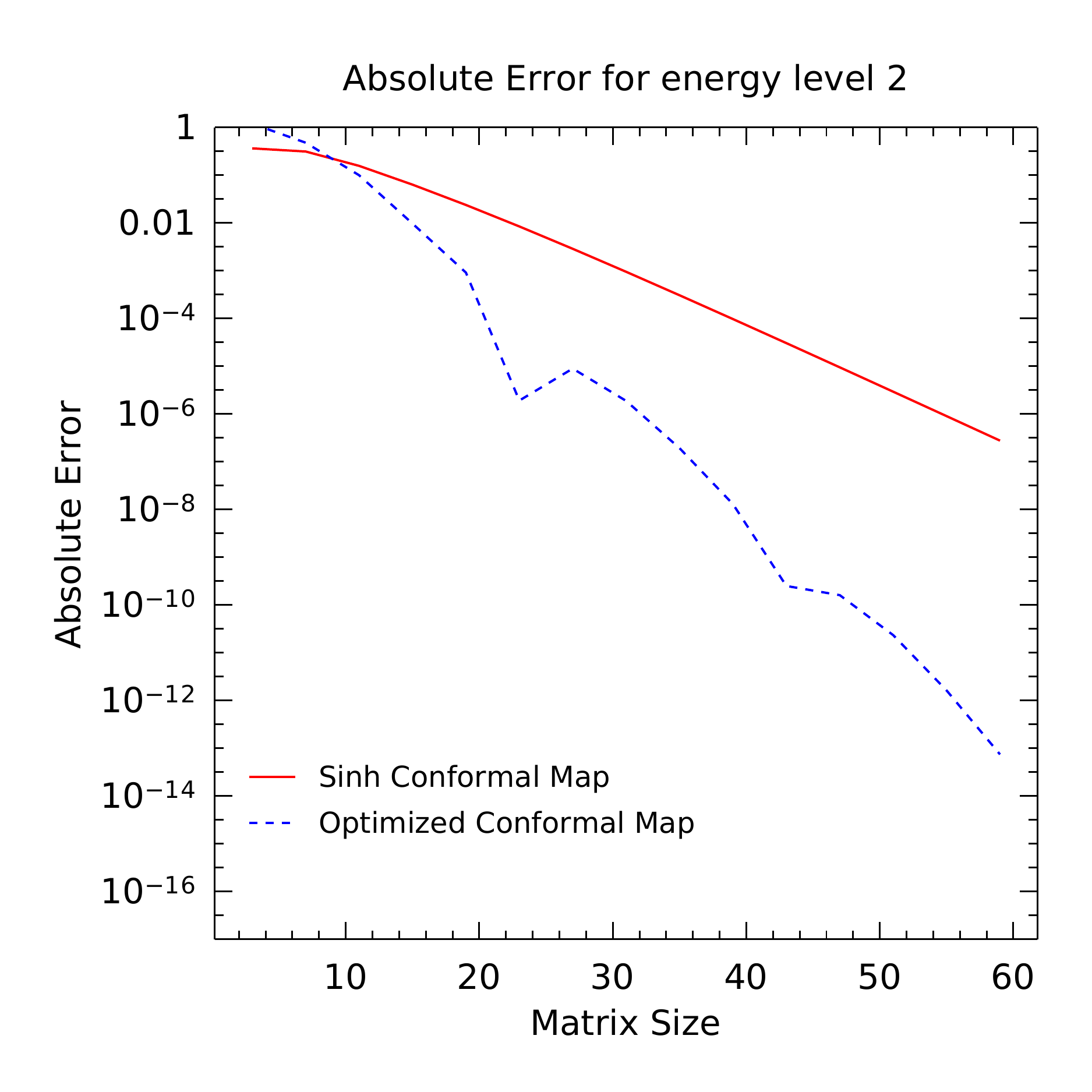} &  \includegraphics[width=0.40\textwidth]{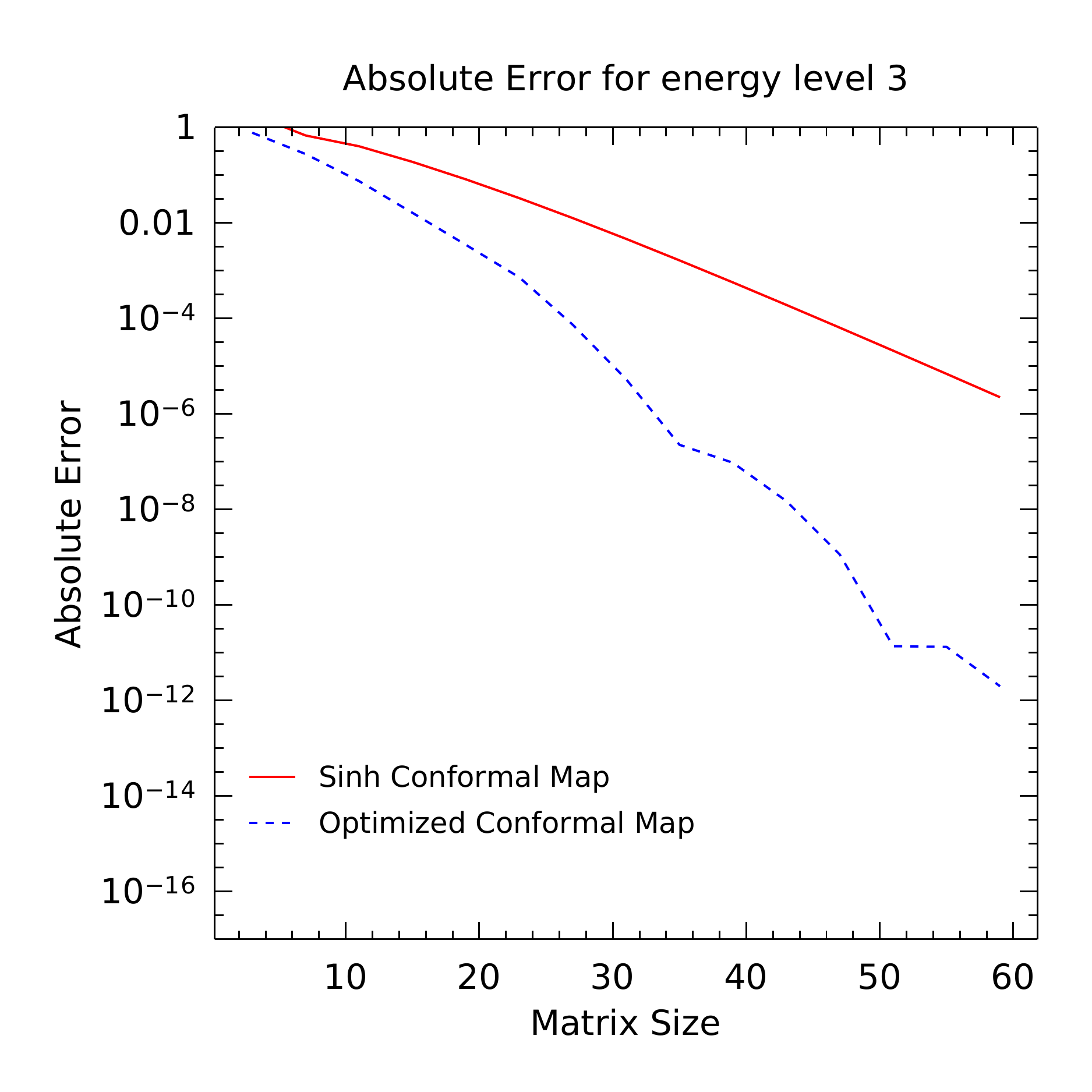} \\
(c) & (d)
\end{tabular}
\caption{Numerical evaluation of the relative error for the potentials $V(x)$ in equation \eqref{formula: exact potential} with $\lambda_{i}(g)$ for $i=1,2,3,4$ as shown in equations \eqref{formula: true value energy} with $\phi(x) = \sinh(x)$. For all figures, we used the value $g=1$. Figure (a) shows the relative error for the potential with $\lambda_{1}(1)  = -6$ with exact eigenvalue $E_{0} = -1$. Figure (b) shows the relative error for the potential with $\lambda_{2}(1)  = -10$ with exact eigenvalue $E_{1} = -3$. Figure (c) shows the relative error for the potential with $\lambda_{3}(1)  =  -13+\sqrt{17}$ with exact eigenvalue $E_{2} = -4+\sqrt{17}$. Figure (d) shows the relative error for the potential with $\lambda_{4}(1)  = -12$ with exact eigenvalue $E_{3} = -1$.}
\label{figure: potentials}
\end{center}
\end{figure}

Using these exact values as comparison, we present the following figures showing the convergence of the DESCM. Figure \ref{figure: potentials}, shows the absolute error between our approximation and the exact values shown in equation \eqref{formula: true value energy}.
In this case, the optimal map can be derived analytically and is given by:
\begin{align}
\phi(t) & = \left(\frac{2}{g}\right)^{1/2}\sinh(t).
\end{align}

Explicitly, we define the absolute error as:
\begin{equation}
{\rm Absolute \, \, error} = \left| {\mathcal E}_{n}(N) - E_{n} \right| \qquad \textrm{for} \quad n =0,1,2,3.
\end{equation}
As we can see from Figure \ref{figure: potentials}, the proposed algorithm converges quite well.

We will now define an approximation to the absolute error as follows:
\begin{equation}\label{formula: approximate relative error}
\epsilon_{n}(N) = \left | {\mathcal E}_{n}(N)-{\mathcal E}_{n}(N-1) \right | \qquad \textrm{for} \quad \left\{\begin{array}{ccc} N &= &1,2,3,\ldots \\
n & = & 0,1,2,\ldots. \end{array} \right.
\end{equation}

We will now consider the general potentials of the form in equation \eqref{formula: rational oscillator}.
\begin{equation}
V(x) = \omega x^{2m} + \dfrac{\displaystyle \sum_{i=0}^{k} \lambda_{i} x^{i}   }{1 + \displaystyle \sum_{j=1}^{2l} g_{j} x^{j}} \qquad \textrm{with} \quad  k-2l<2m,
\end{equation}
These potentials contains many free parameters. Indeed, we have $5+k+2l$ free parameters including:
\begin{enumerate}
\item The exponent $m$
\item The coefficient parameters $\omega$, $\lambda_{i}$ and $g_{j}$ for $i=0,\ldots,k$ and $j=1,\ldots,2l$
\item The degrees $k$ and $2\,l$ of the polynomials in the numerator and denominator.
\end{enumerate}

It is important to mention that choosing the parameters $g_{i}, i=1,\ldots,2l$ at random such that $q(x) \neq 0,  \forall x\in\mathbb{R}$ might at first seem problematic. However, by the fundamental theorem of algebra, we are assured that $q(x)$ may be written as:
\begin{equation}\label{q(x) polynomial}
q(z) = 1 + \displaystyle \sum_{i=1}^{2l} g_{i} z^{i} = \dfrac{\displaystyle\prod_{i=1}^{l} (z-z_{i})(z-\bar{z_{i}})}{\displaystyle\prod_{i=1}^{l} z_{i}\bar{z_{i}}},
\end{equation}
where $\{z_{i},\bar{z_{i}}\}_{i=1}^{l}$ are the complex zeros of the strictly positive polynomial $q(x)$. In this cases, $\bar{z_{i}}$ denotes the complex conjugate of $z_{i}$. Hence, the polynomials $q(z)$ can be easily constructed by choosing values for $\{\Re\{z_{i}\},\Im\{z_{i}\}\}_{i=1}^{l}$ and substituting them into equation \eqref{q(x) polynomial}.

To remove as much bias as possible in the calculation, we let the parameters involved in the computation represent random variables such that:
\begin{align}\label{random variables}
\omega & \sim U(0,10) \nonumber\\
k & \sim U\{0,1,\ldots,2m+2l-1\} \nonumber \\
\lambda_{i} & \sim U(-10,10) \qquad i =0,\ldots,k \\
\Re\{z_{i}\} & \sim U(-5,5) \qquad i =1,\ldots,l \nonumber \\
\Im\{z_{i}\} & \sim U(0,10) \qquad i =1,\ldots,l \nonumber,
\end{align}
for fixed parameters $m$ and $l$ and where $\Re \left\{z\right\}$ stands for the real part of $z$.

We acknowledge that we have already used the symbol $"\sim"$ to denote the concept of asymptoticity. However, in equation \eqref{random variables}, the symbol $X\sim f$ denotes the common statistics notation that the random variable $X$ follows the distribution $f$.  In \eqref{random variables}, $U(a,b)$ denotes the continuous uniform distribution on the interval $(a,b)$ and $U\{c,c+1,\ldots,d\}$ denotes the discrete uniform distribution for the integer support $\{c,c+1,\ldots,d\}$, $c,d\in\mathbb{N}_{0}$.

\begin{figure}[!h]
\begin{center}
\begin{tabular}{cccc} \includegraphics[width=0.50\textwidth]{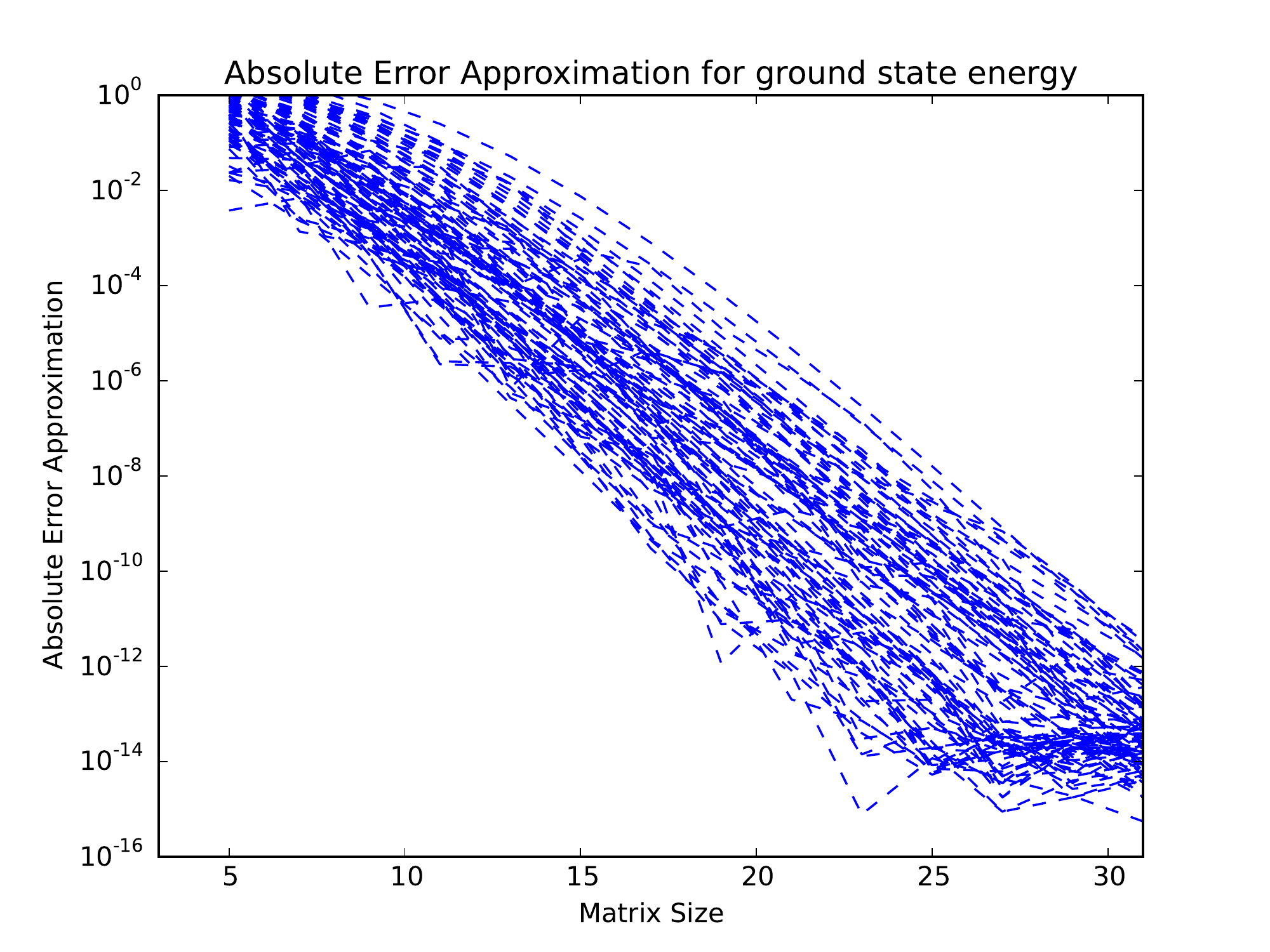} &  \includegraphics[width=0.50\textwidth]{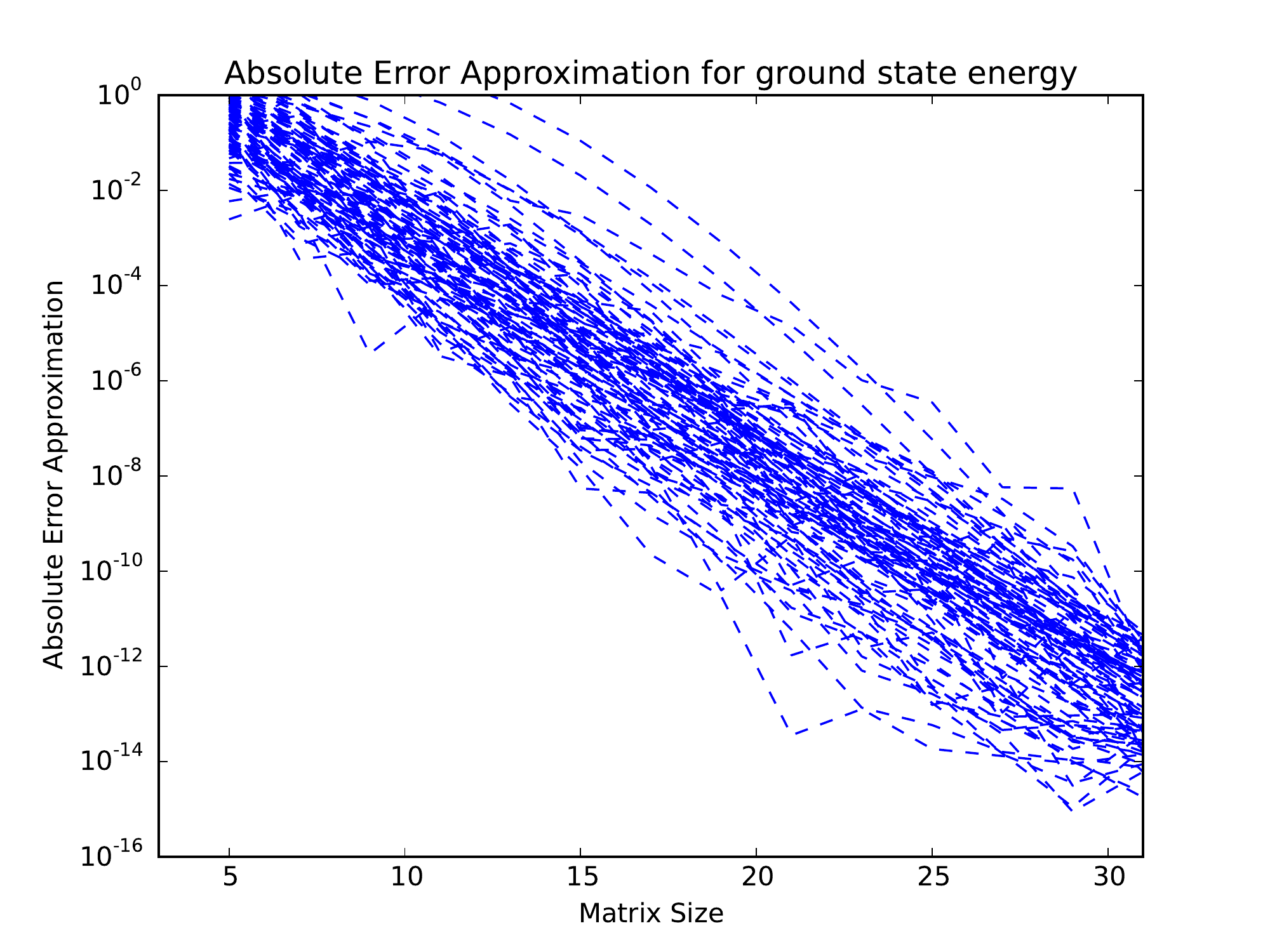} \\
(a) $m=1$ & (b) $m=2$ \\
\includegraphics[width=0.50\textwidth]{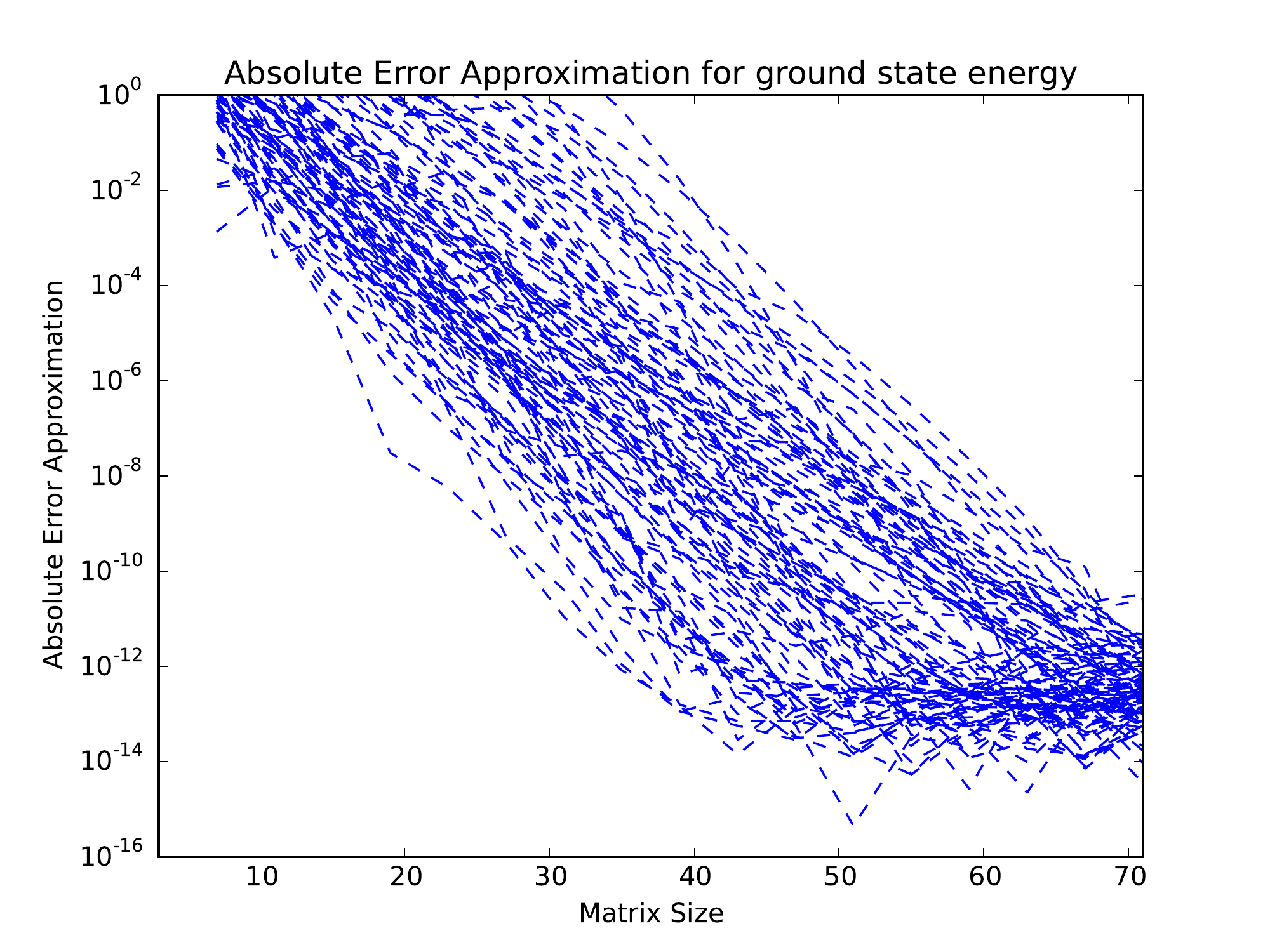} &  \includegraphics[width=0.50\textwidth]{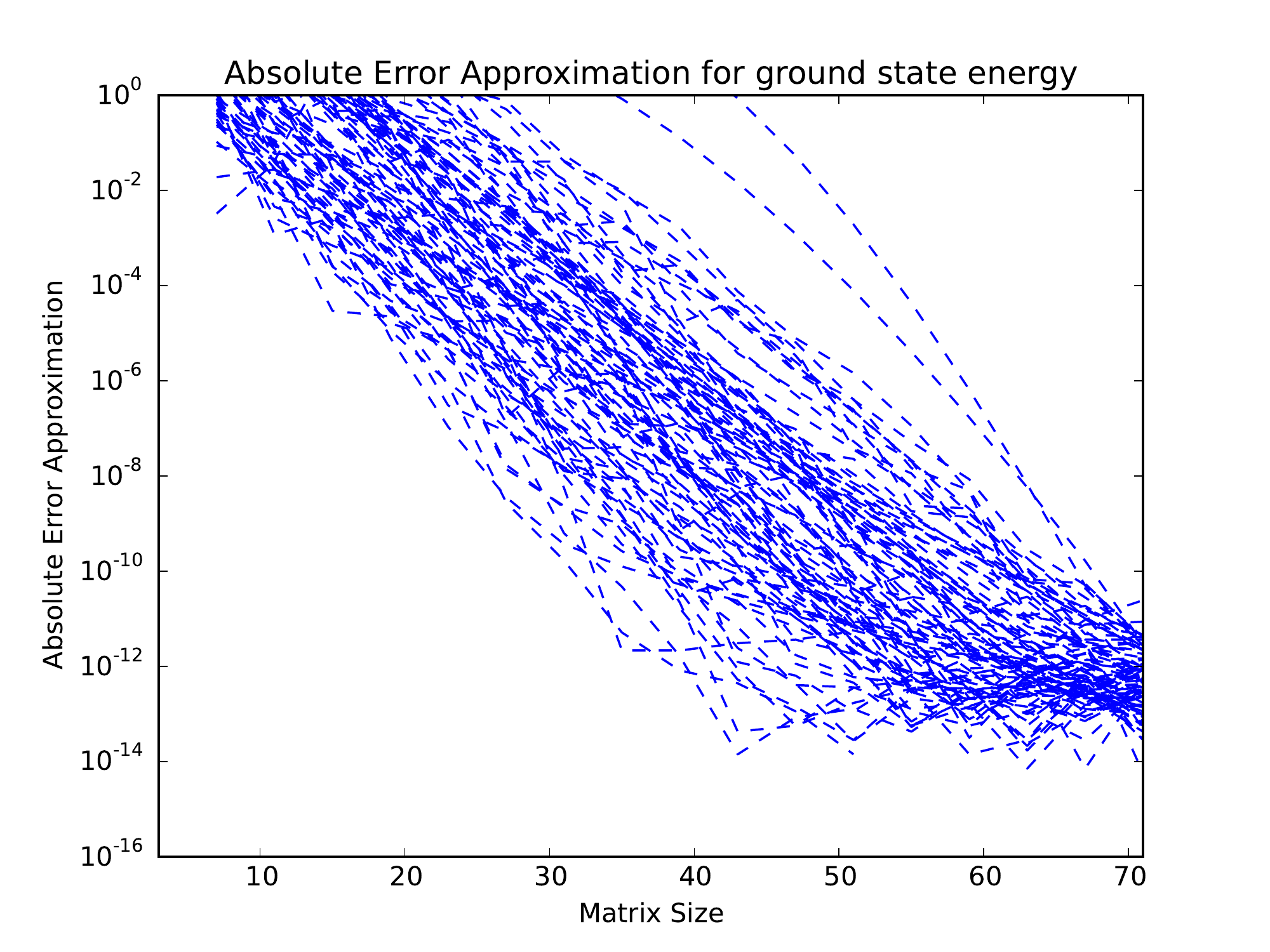}
 \\
(c) $m=3$ & (d) $m=4$
\end{tabular}
\caption{Application of the DESCM for 100 randomly generated potentials of the form \newline
 $V(x) =\omega x^{2m} + \dfrac{ \sum_{i=0}^{k} \lambda_{i} x^{i}   }{1 +  g_{1} x+ g_{2} x^2} $. Figure (a) corresponds to $m=1$. Figure (b) corresponds to $m=2$. Figure (c) corresponds to $m=3$. Figure (d) corresponds to $m=4$.}
\label{figure: potentials1}
\end{center}
\end{figure}

In Figure \ref{figure: potentials1}, we applied the DESCM to 100 randomly generated potentials  of the form in \eqref{formula: rational oscillator} according to \eqref{random variables} with $m=1,2,3,4$ and $l=1$. In this case, we can also find the optimal map analytically which is given by:
\begin{align}
\phi(t) & = \left[ \Im\{z_{1}\}\csc\left( \frac{\pi}{4}\right) \right]\sinh(t) + \Re\{z_{1}\}.
\end{align}
As we can, the DESCM performs quite well for a wide variety of parameters values.

\begin{figure}[!h]
\begin{center}
\begin{tabular}{cccc} \includegraphics[width=0.50\textwidth]{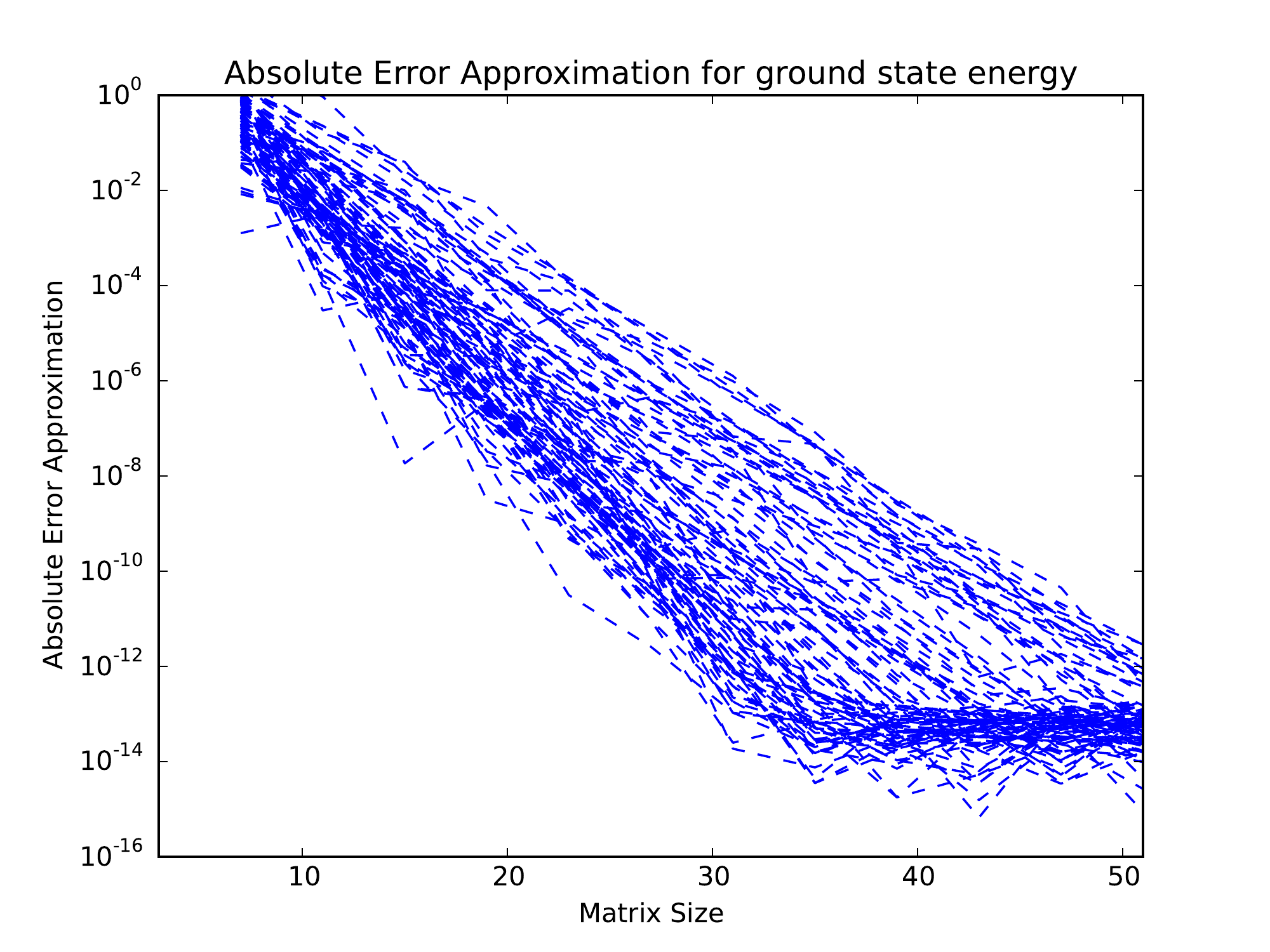} &  \includegraphics[width=0.50\textwidth]{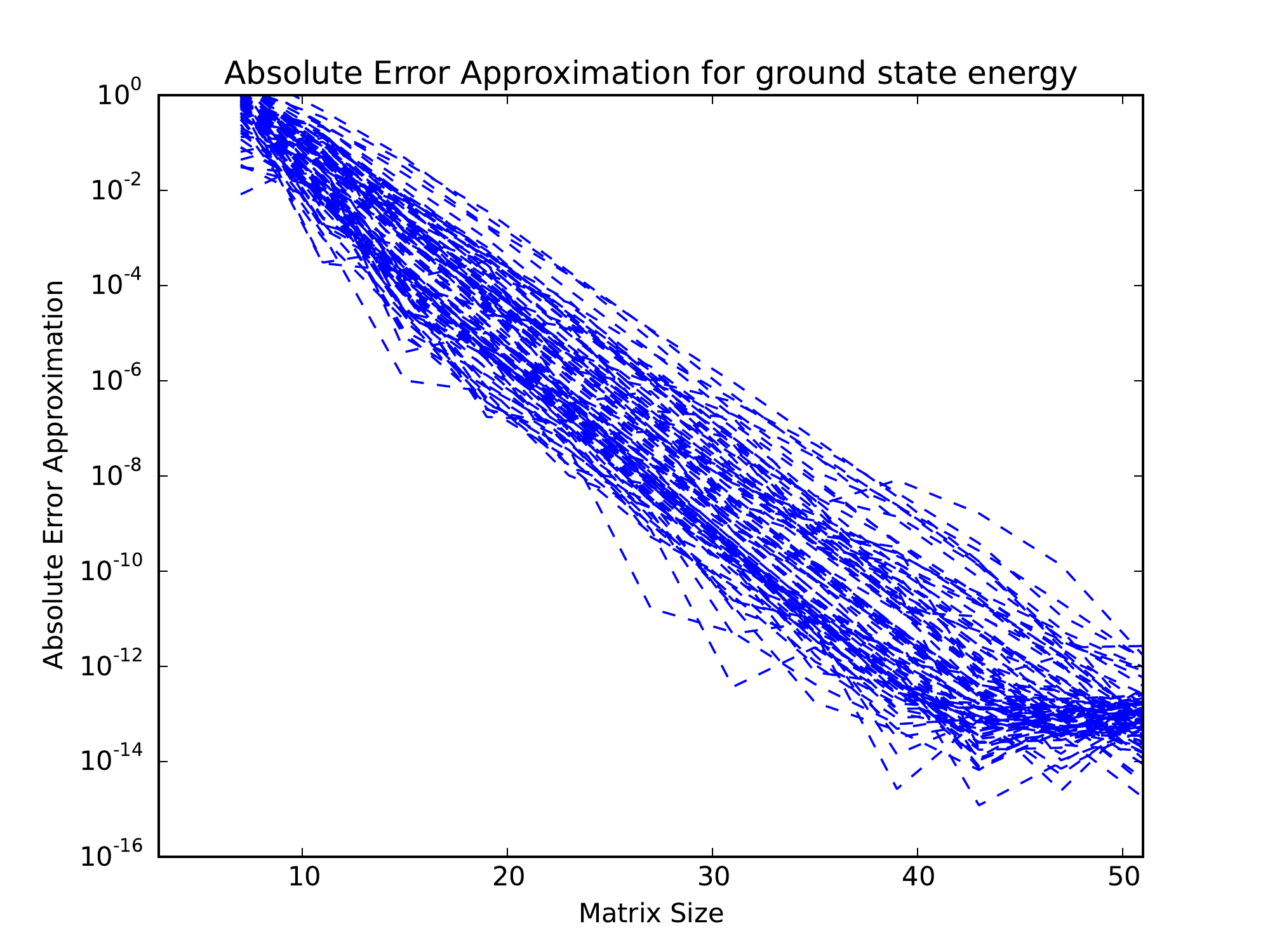}\\
(a) $m=1$ & (b) $m=2$ \\
\includegraphics[width=0.50\textwidth]{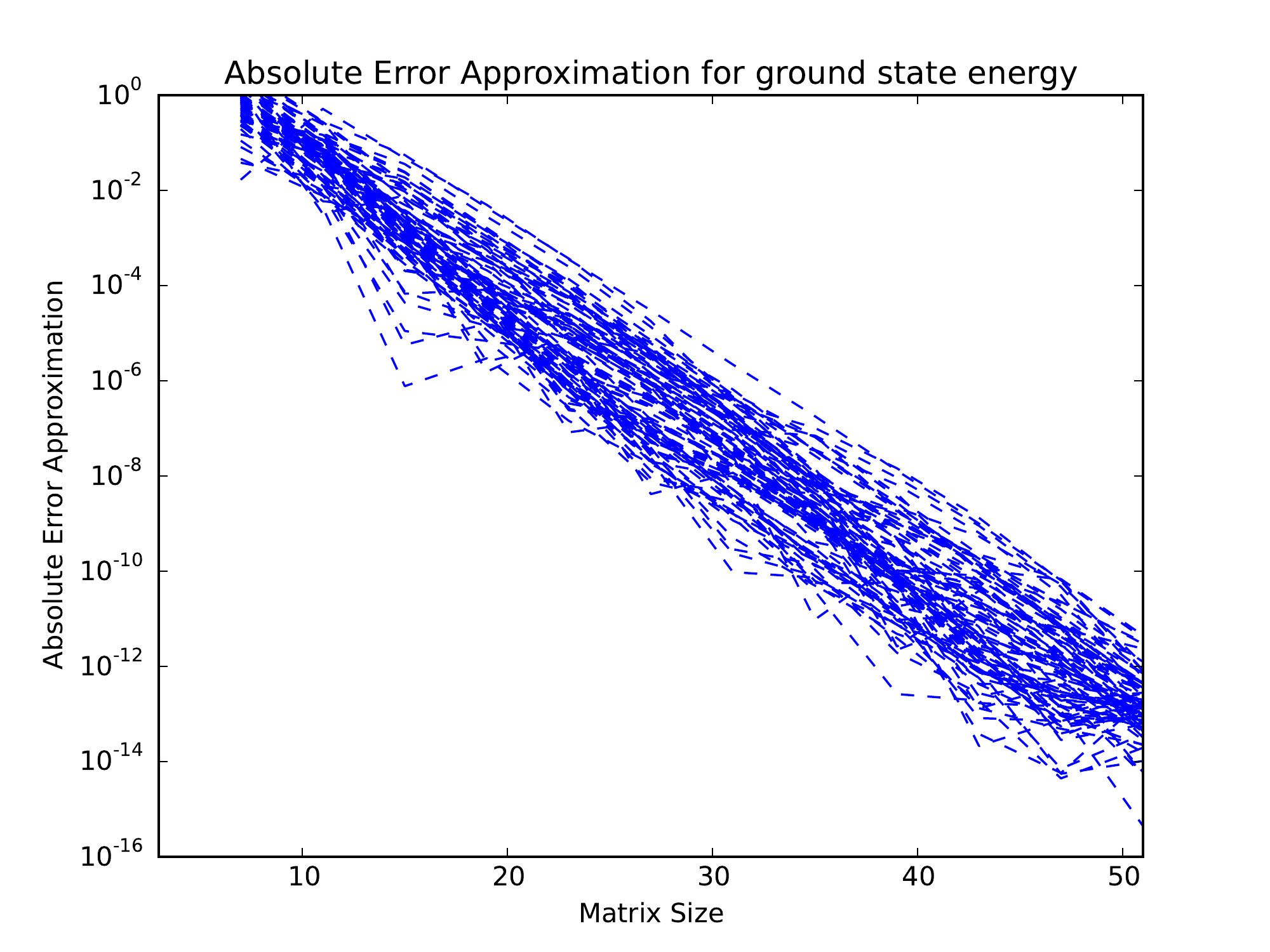} &  \includegraphics[width=0.50\textwidth]{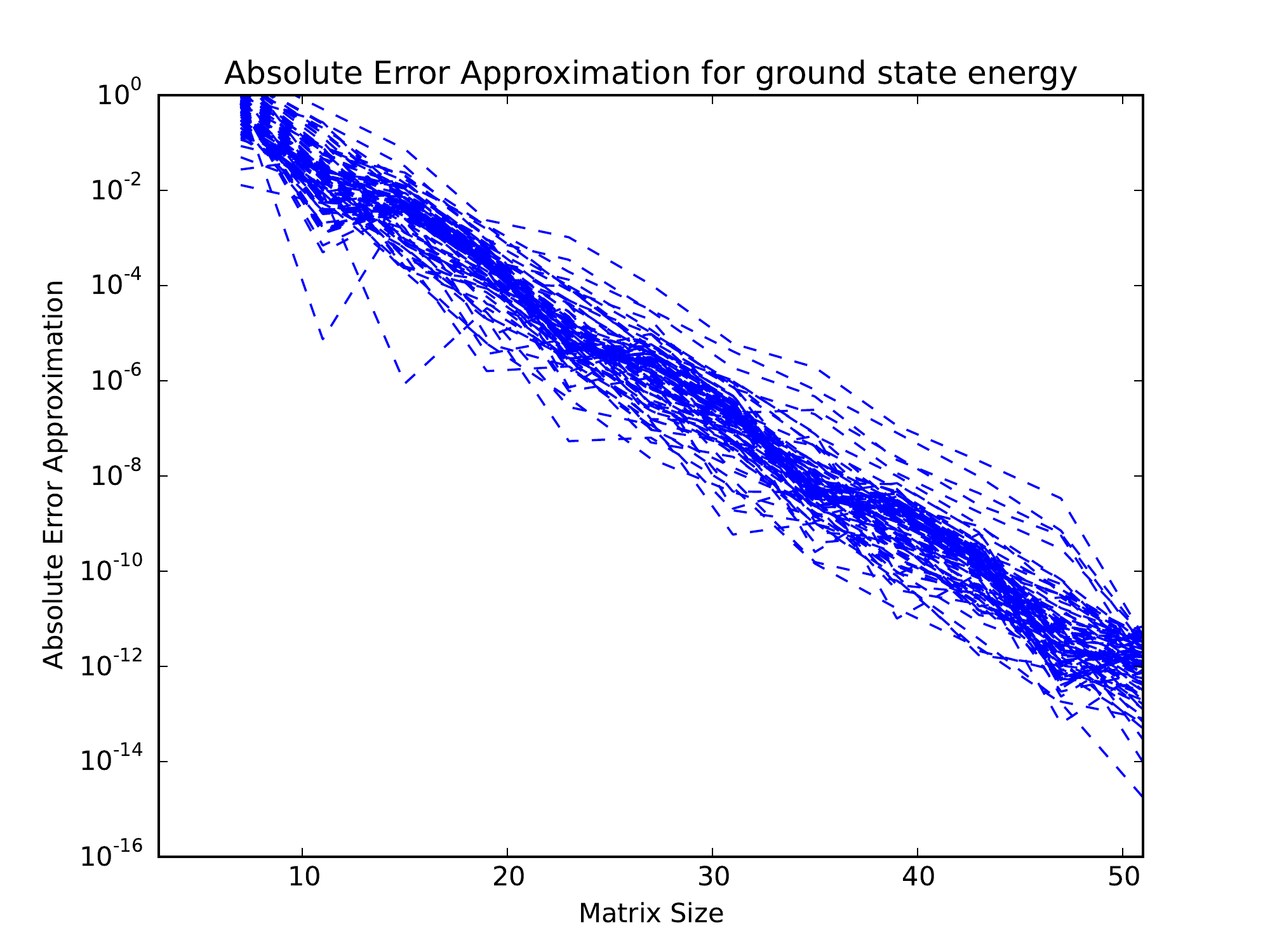}
 \\
(c) $m=3$ & (d) $m=4$
\end{tabular}
\caption{Application of the DESCM for 100 randomly generated potentials of the form \newline
 $V(x) =\omega x^{2m} + \dfrac{ \sum_{i=0}^{k} \lambda_{i} x^{i}   }{1 +  g_{1} x+ g_{2} x^2+g_{3} x^3+g_{4} x^4} $. Figure (a) corresponds to $m=1$. Figure (b) corresponds to $m=2$. Figure (c) corresponds to $m=3$. Figure (d) corresponds to $m=4$.}
\label{figure: potentials2}
\end{center}
\end{figure}

In Figure \ref{figure: potentials2}, we applied the DESCM to 100 randomly generated potentials  of the form in \eqref{formula: rational oscillator} according to \eqref{random variables} with $m=1,2,3,4$ and $l=2$. Unlike the previous examples, the optimal map cannot be found analytically for these types of potentials. This would lead us to solve the non-linear program in \eqref{formula: system of equations} for every randomly generated potential. However, this is much harder than anticipated because this non-linear program must be calibrated for any potential. More explicitly, solving the non-linear program requires a user input of two parameters, "obj$\_$scaling$\_$factor" which controls the scaling of the objective function and "Hint" which controls the homotopy solution process for the nonlinear program. Given the sensitivity of non-linear programming, these parameters must be tuned for any arbitrary potential. As such, it is an infeasible task to finely tune these parameters for any given list of randomly generated potentials. However, to demonstrate the power of the DESCM, Figure \ref{figure: potentials2} display the implementation using the basic non-optimal $\phi(t) = \sinh(t)$ conformal mapping. As we see, the DESCM still converges quite well for a wide range of parameters.

\section{Conclusion}
Several methods have been used to evaluate the energy eigenvalues of perturbed harmonic oscillators by rational potentials. In this work, we present a method based on the DESCM where the wave function of a transformed Schr\"odinger equation \eqref{formula: transformed Schrodinger equation} is approximated by as a Sinc expansion. By summing over $2N+1$ collocation points, the implementation of the DESCM leads to a generalized eigenvalue problem with symmetric and positive definite matrices. In addition, we also show that the convergence of the DESCM in this case can be improved to the rate ${\cal O} \left( \left(\dfrac{N^{5/2}}{\log(N)^{2}} \right)  \exp \left(- \kappa^{\prime}\dfrac{N}{\log(N)} \right)\right)$ as $N\to \infty$ where $2N+1$ is the dimension of the resulting generalized eigenvalue system and $\kappa^\prime$ is a constant that depends on the potential. The convergence of this method can be improved by adding a polynomial adjustment to the typical $\sinh$ conformal mapping to displace the complex singularities away from the real axis.

\section{Tables and Figures}
\clearpage

\clearpage
%\bibliography{C:/Users/Philippe/Documents/Research/MyBibliography/library}
%\bibliography{/Users/hsafouhi/Documents/DATA1/01 AMybibliography/library}

%\bibliography{library}

\begin{thebibliography}{10}

\bibitem{Mitra1978}
A.K. Mitra.
\newblock {On the interaction of the type $\lambda x^2/(1+g x^2)$}.
\newblock {\em Journal of Mathematical Physics}, 19(10):2018--2022, 1978.

\bibitem{Summers2012}
S.J. Summers.
\newblock {A Perspective on Constructive Quantum Field Theory}.
\newblock {\em arXiv:1203.3991}, page~59, 2012.

\bibitem{Magnen2008}
J.~Magnen and V.~Rivasseau.
\newblock {Constructive $\phi^4$ field theory without tears}.
\newblock {\em Annales Henri Poincare}, 9(2):403--424, 2008.

\bibitem{Rivasseau2012a}
V.~Rivasseau and Z.~Wang.
\newblock {Constructive renormalization for $\Phi^{4}_{2}$ theory with loop
  vertex expansion}.
\newblock {\em Journal of mathematical physics}, 53(4):042302--1 -- 042302--29,
  2012.

\bibitem{Jelic2012}
V.~Jelic and F.~Marsiglio.
\newblock {The double-well potential in quantum mechanics: a simple,
  numerically exact formulation}, 2012.

\bibitem{Griffiths2004}
D.J. Griffiths.
\newblock {\em {Introduction to Quantum Mechanics}}.
\newblock Addison-Wesley, 2nd editio edition, 2004.

\bibitem{Anderlini2007}
M.~Anderlini, P.J. Lee, B.L. Brown, J.~Sebby-Strabley, W.D. Phillips, and J.V.
  Porto.
\newblock {Controlled exchange interaction between pairs of neutral atoms in an
  optical lattice.}
\newblock {\em Nature}, 448(7152):452--456, 2007.

\bibitem{Trotzky2008}
S.~Trotzky, P.~Cheinet, S.~F{\"{o}}lling, M.~Feld, U.~Schnorrberger, A.M. Rey,
  A.~Polkovnikov, E.A. Demler, M.D. Lukin, and I.~Bloch.
\newblock {Time-resolved observation and control of superexchange interactions
  with ultracold atoms in optical lattices.}
\newblock {\em Science (New York, N.Y.)}, 319(5861):295--299, 2008.

\bibitem{Foot2009}
C.J. Foot and M.~Shotter.
\newblock {Double well potentials and quantum gates}.
\newblock {\em American Journal of Physics}, 79(7):762--768, 2011.

\bibitem{Murmann2014}
S.~Murmann, A.~Bergschneider, V.M. Klinkhamer, G.~Z{\"{u}}rn, T.~Lompe, and
  S.~Jochim.
\newblock {Two Fermions in a double well : Exploring a fundamental building
  block of the Hubbard model}.
\newblock {\em arXiv:1410.8784}, 2014.

\bibitem{Bessis1980}
N.~Bessis and G.~Bessis.
\newblock {A note on the Schr{\"o}dinger equation for the $x^2+\lambda x^2/(1+g
  x^2)$ potential}.
\newblock {\em Journal of Mathematical Physics}, 21(12):2780--2785, 1980.

\bibitem{Bessis1983}
N.~Bessis, G.~Bessis, and G.~Hadinger.
\newblock {Perturbed harmonic oscillator ladder operators: eigenenergies and
  eigenfunctions for the $x^2 + \lambda x^2 /(1+g x^2)$ interaction}.
\newblock {\em Journal of Physics A: Mathematical and General}, 16(3):497--512,
  1983.

\bibitem{Bessis1992}
N.~Bessis and G.~Bessis.
\newblock {Perturbed factorization of the symmetric-anharmonic-oscillator
  eigenequation}.
\newblock {\em Physical Review A}, 46(11):6824--6845, 1992.

\bibitem{Fack1986}
V.~Fack, H.E. {De Meyer}, and G.V. Berghe.
\newblock {Dynamic-group approach to the $x^2+\lambda x^2/(1+g x^2)$
  potential}.
\newblock {\em Journal of Mathematical Physics}, 27(5):1340, 1986.

\bibitem{Stubbins1995}
C.~Stubbins and M.~Gornstein.
\newblock {Variational estimates of the energies for the potential $x^2 +
  \lambda x^2 /(1+g x^2)$}.
\newblock {\em Physics Letters A}, 202(1):34--39, 1995.

\bibitem{Chaudhuri1983}
R.N. Chaudhuri and B.~Mukherjee.
\newblock {On the Schr{\"o}dinger equation for the interaction $x^2 + \lambda
  x^2 /(1+g x^2 )$}.
\newblock {\em Journal of Physics A: Mathematical and General},
  16(17):4031--4038, 1983.

\bibitem{Blecher1987}
M.H. Blecher and P.G.L. Leach.
\newblock {The Schr{\"o}dinger equation for the $x^2 +\lambda x^2 /(1+g x^2)$
  interaction}.
\newblock {\em Journal of Physics A: Mathematical and General},
  20(17):5923--5927, 1987.

\bibitem{Filho1989}
E.D. Filho and R.M. Ricotta.
\newblock {Supersymmetric quantum mechanics and higher excited states of a
  non-polynomial potential}.
\newblock {\em Modern Physics Letters A}, 4(23):2283--2288, 1989.

\bibitem{Filho1991}
E.D. Filho and R.M. Ricotta.
\newblock {Partial algebraization of the non-polynomial potential}.
\newblock {\em Modern Physics Letters A}, 6(23):2137--2142, 1991.

\bibitem{Flessas1982a}
G.P. Flessas.
\newblock {Definite integrals as solutions for the $x^2 +\lambda x^2 /(1+ g
  x^2)$ potential}.
\newblock {\em Journal of Physics A: Mathematical and General},
  15(3):L97--L101, 1982.

\bibitem{Flessas1981}
G.P. Flessas.
\newblock {On the Schr{\"o}dinger equation for the $x^2 +\lambda x^2 /(1+ g
  x^2)$ interaction}.
\newblock {\em Physics Letters A}, 83(3):121--122, 1981.

\bibitem{Heading1982}
J.~Heading.
\newblock {Polynomial-type eigenfunctions}.
\newblock {\em Journal of Physics A: Mathematical and General},
  15(8):2355--2367, 1982.

\bibitem{Gallas1988}
J.A.C. Gallas.
\newblock {Exact analytical eigenfunctions for the $x^2 + \lambda x^2 /(1+g x^2
  )$ interaction}.
\newblock {\em Journal of Physics A: Mathematical and General},
  21(16):3393--3397, 1988.

\bibitem{Marcilhacy1985}
G.~Marcilhacy and R.~Pons.
\newblock {The Schr{\"o}dinger equation for the interaction potential $x^2
  +\lambda x^2 /(1+g x^2 )$ and the first Heun confluent equation}.
\newblock {\em Journal of Physics A: Mathematical and General},
  18(13):2441--2449, 1985.

\bibitem{Roy1987}
P.~Roy and R.~Roychoudhury.
\newblock {New integral solutions of the non-polynomial oscillator $V(x)=x^2
  +\lambda x^2 /(1+g x^2)$ when $\lambda=2g(2-3g)$}.
\newblock {\em Journal of Physics A: Mathematical and General},
  20(18):L1245--L1248, 1987.

\bibitem{Roy1987a}
P.~Roy and R.~Roychoudhury.
\newblock {New exact solutions of the non-polynomial oscillator $x^2 +\lambda
  x^2 /(1+ g x^2)$ and supersymmetry}.
\newblock {\em Physics Letters A}, 122(6-7):275--279, 1987.

\bibitem{Roy1988}
P.~Roy, R.~Roychoudhury, and Y.P. Varshni.
\newblock {Some solutions of a supersymmetric nonpolynomial oscillator-a
  comparison between the SWKB and WKB methods}.
\newblock {\em Journal of Physics A: Mathematical and General},
  21(7):1589--1594, 1988.

\bibitem{Roy1990}
P.~Roy and R.~Roychoudhury.
\newblock {Some observations on the nature of solutions for the interaction
  $V(x)=x^2 +(\lambda x^2 /(1+gx 2 ))$}.
\newblock {\em Journal of Physics A: Mathematical and General},
  23(9):1657--1660, 1990.

\bibitem{Berghe1989}
G.V. Berghe and H.E. {De Meyer}.
\newblock {Pairs of analytical eigenfunctions for the $x^2 + \lambda x^2 /(1 +
  gx^2)$ interaction}.
\newblock {\em Journal of Physics A: Mathematical and General},
  22(10):1705--1710, 1989.

\bibitem{Varma1981}
V.S. Varma.
\newblock {On the $x^2 + \lambda x^2 /(1+g x^2)$ interaction}.
\newblock {\em Journal of Physics A: Mathematical and General},
  14(12):L489--L492, 1981.

\bibitem{Whitehead1982}
R.R. Whitehead, A.~Watt, G.P. Flessas, and M.A. Nagarajan.
\newblock {Exact solutions of the Schr{\"o}dinger equation $\left(
  -\frac{d^{2}}{d x^2} + x^2 + \frac{\lambda x^2}{1+g x^2} \right) \psi(x) = E
  \psi(x) $}.
\newblock {\em Journal of Physics A: Mathematical and General},
  15(4):1217--1226, 1982.

\bibitem{Estrin1990a}
D.A. Estr{\'{\i}}n, F.M. Fern{\'{a}}ndez, and E.A. Castro.
\newblock {On the Hill determinant method}.
\newblock {\em Journal of Physics A: Mathematical and General},
  23(12):2395--2400, 1990.

\bibitem{Agrawal1993}
R.K. Agrawal and V.S. Varma.
\newblock {Rational potential using a modified Hill determinant method}.
\newblock {\em Physical Review A}, 48(3):1921--1928, 1993.

\bibitem{Hautot1981}
A.~Hautot.
\newblock {Interaction $\lambda x^2 (1+ g x^2)$ revisited}.
\newblock {\em Journal of Computational Physics}, 39(1):72--93, 1981.

\bibitem{Fernandez1991}
F.M. Fern{\'{a}}ndez.
\newblock {Convergent power-series solutions to the Schr{\"o}dinger equation
  with the potential}.
\newblock {\em Physics Letters A}, 160(2):116--118, 1991.

\bibitem{Hodgson1988}
R.J.W. Hodgson.
\newblock {High-precision calculation of the eigenvalues for the $x^2 +\lambda
  x^2 /(1+g x^2)$ potential}.
\newblock {\em Journal of Physics A: Mathematical and General},
  21(7):1563--1570, 1988.

\bibitem{Fassari1996}
S.~Fassari.
\newblock {A note on the eigenvalues of the Hamiltonian of the harmonic
  oscillator perturbed by the potential $\lambda x^2/(1+g x^2)$}.
\newblock {\em Reports on Mathematical Physics}, 37(2):283--293, 1996.

\bibitem{Fassari1997}
S.~Fassari and G.~Inglese.
\newblock {On the eigenvalues of the Hamiltonian of the harmonic oscillator
  with the interaction lambda $x^2 / (1+g x^2)$ (II)}.
\newblock {\em Reports on Mathematical Physics}, 39(1):77--86, 1997.

\bibitem{Kaushal1979}
R.S. Kaushal.
\newblock {Small $g$ and large $\lambda$ solution of the Schr{\"o}dinger
  equation for the interaction $\lambda x^2 /(1+g x^2 )$}.
\newblock {\em Journal of Physics A: Mathematical and General},
  12(10):L253--L258, 1979.

\bibitem{Lai1982}
C.S. Lai and H.E. Lin.
\newblock {On the Schr{\"o}dinger equation for the interaction $x^2 +\lambda
  x^2 /(1+g x^2 )$}.
\newblock {\em Journal of Physics A: Mathematical and General},
  15(5):1495--1502, 1982.

\bibitem{Fack1987}
V.~Fack, G.V. Berghe, and H.E. {De Meyer}.
\newblock {Some finite difference methods for computing eigenvalues and
  eigenvectors of special two-point boundary value problems}.
\newblock {\em Journal of Computational and Applied Mathematics},
  20(1):211--217, 1987.

\bibitem{Galicia1979}
S.~Galicia and J.~Killingbeck.
\newblock {Accurate calculation of perturbed oscillator energies}.
\newblock {\em Physics Letters A}, 71(1):17--18, 1979.

\bibitem{Ishikawa2007}
H.~Ishikawa.
\newblock {Numerical methods for the eigenvalue determination of second-order
  ordinary differential equations}.
\newblock {\em Journal of Computational and Applied Mathematics},
  208(2):404--424, 2007.

\bibitem{Bhagwat1981}
K.V. Bhagwat.
\newblock {A harmonic oscillator perturbed by the potential $\lambda x^2 /(1+g
  x^2 )$}.
\newblock {\em Journal of Physics A: Mathematical and General}, 14(2):377--378,
  1981.

\bibitem{Simos1998}
T.E. Simos.
\newblock {Some embedded modified Runge-Kutta methods for the numerical
  solution of some specific Schr{\"o}dinger equations}.
\newblock {\em Journal of mathematical chemistry}, 24(1-3):23--37, 1998.

\bibitem{Simos1997}
T.E. Simos.
\newblock {A new finite-difference scheme for the numerical solution of the
  Schr{\"o}dinger equation}.
\newblock {\em Canadian Journal of Physics}, 75(5):325--335, 1997.

\bibitem{Ying2010}
H.~Ying, Z.~Fan-Ming, Y.~Yan-Fang, and L.~Chun-Fang.
\newblock {Energy eigenvalues from an analytical transfer matrix method}.
\newblock {\em Chinese Physics B}, 19(4):040306--1 -- 040306--6, 2010.

\bibitem{DaCosta2008}
G.A.T.F. da~Costa and M.~Gomes.
\newblock {Borel-Leroy summability of a nonpolynomial potential}.
\newblock {\em Reports on Mathematical Physics}, 61(3):401--415, 2008.

\bibitem{Hislop1990}
D.~Hislop, M.F. Wolfaardt, and P.G.L. Leach.
\newblock {The Schr{\"o}dinger equation for the $f(x)/g(x)$ interaction}.
\newblock {\em Journal of Physics A: Mathematical and General},
  23(21):L1109--L1112, 1990.

\bibitem{Nanayakkara2012}
A.~Nanayakkara.
\newblock {Asymptotic behavior of eigenenergies of nonpolynomial oscillator
  potentials $V (x) = x^{2N} + ( \lambda x^{m_{1}} )/(1 + g x^{m_{2}} )$}.
\newblock {\em Canadian Journal of Physics}, 90(6):585--592, 2012.

\bibitem{Scherrer1988}
H.~Scherrer, H.~Risken, and T.~Leiber.
\newblock {Eigenvalues of the Schr{\"o}dinger equation with rational
  potentials}.
\newblock {\em Physical Review A}, 38(8):3949--3959, 1988.

\bibitem{Witwit1996a}
M.R.M. Witwit.
\newblock {Energy levels of a nonpolynomial oscillator using finite difference
  technique}.
\newblock {\em Journal of Computational and Applied Mathematics},
  69(2):331--343, 1996.

\bibitem{Gaudreau2014a}
P.~Gaudreau, R.M. Slevinsky, and H.~Safouhi.
\newblock {The Double Exponential Sinc Collocation Method for Singular
  Sturm-Liouville Problems}.
\newblock {\em arXiv:1409.7471v2}, 2014.

\bibitem{Stenger1981}
F.~Stenger.
\newblock {Numerical Methods Based on Whittaker Cardinal, or Sinc Functions}.
\newblock {\em SIAM Review}, 23(2):165--224, 1981.

\bibitem{Gaudreau2015b}
P.~Gaudreau, R.M. Slevinsky, and H.~Safouhi.
\newblock {Computing energy eigenvalues of anharmonic oscillators using the
  double exponential Sinc collocation method}.
\newblock {\em Annals of Physics}, 360:520--538, 2015.

\bibitem{Eggert1987}
N.~Eggert, M.~Jarratt, and J.~Lund.
\newblock {Sinc function computation of the eigenvalues of Sturm-Liouville
  problems}.
\newblock {\em Journal of Computational Physics}, 69(1):209--229, 1987.

\bibitem{Stenger1979}
F.~Stenger.
\newblock {A "Sinc-Galerkin" method of solution of boundary value problems}.
\newblock {\em Mathematics of Computation}, 33(145):85--109, 1979.

\bibitem{Slevinsky2014a}
R.M. Slevinsky and S.~Olver.
\newblock {On the use of conformal maps for the acceleration of convergence of
  the trapezoidal rule and Sinc numerical methods}.
\newblock arXiv(1406.3320):1--25, 2014.

\bibitem{Bezanson2012}
J.~Bezanson, S.~Karpinski, V.B. Shah, and A.~Edelman.
\newblock {Julia: A Fast Dynamic Language for Technical Computing}.
\newblock arXiv(1209.5145):1--27, 2012.

\bibitem{Anderson1999}
E.~Anderson, Z.~Bai, C.~Bischof, S.~Blackford, J.~Demmel, J.~Dongarra, J.~{Du
  Croz}, A.~Greenbaum, S.~Hammarling, A.~McKenney, and D.~Sorensen.
\newblock {\em {{LAPACK} Users' Guide}}.
\newblock Society for Industrial and Applied Mathematics, Philadelphia, PA,
  third edition, 1999.

\end{thebibliography}

\end{document}